\documentclass[12pt]{article}

\usepackage{enumerate}
\usepackage{amsmath}
\usepackage{amssymb,latexsym}
\usepackage{amsthm}

\usepackage{anysize}
\usepackage{graphicx}
\usepackage{url}
\usepackage{color}
\usepackage[english]{babel}
\usepackage{latexsym}
\usepackage{amssymb,amsthm,amsmath}
\usepackage{graphicx,color}
\usepackage{verbatim,setspace}
\newtheorem{theorem}{Theorem}[section]
\newtheorem{proposition}[theorem]{Proposition}
\newtheorem{lemma}[theorem]{Lemma}

\newtheorem{remark}[theorem]{Remark}
\newcommand{\Aut}{{\rm Aut}}
\newcommand{\Paut}{{\rm Paut}}

\newcommand{\diag}{{\rm diag}}
\newcommand{\ord}{{\rm ord}}

\newcommand{\Hess}{{\rm Hess}}
\newcommand{\Ima}{{\rm Im}}

\newcommand{\rank}{{\rm rank }}

\newcommand{\Ker}[1]{\mbox{${\rm Ker\ }{#1}$}}

\newcommand{\AAA}[1]{\mbox{${\bf A}_{#1}$}}
\newcommand{\SSS}[1]{\mbox{${\bf S}_{#1}$}}

\newcommand{\Proof}{{\it Proof}.\ }

\begin{document}

\title{$\SSS{5}$-invariant Nonsingular Quartic Surfaces
%\thanks{This research was  carried out with
%the support of the Italian MIUR (progetto "Strutture
%Geometriche, Algebriche e Combinatoria"), and of
%GNSAGA.}
}
%\titlerunning{Characterization of the Fermat curve}

\author{
Giorgio Faina, Stefano Marcugini and Fernanda Pambianco
\thanks{The research was supported by the Italian Ministero dell'Istruzione, dell'Universit\`a e della Ricerca (MIUR) and
by the Gruppo Nazionale per le Strutture Algebriche, Geometriche e le loro Applicazioni (GNSAGA). The principal results of this paper have been obtained during H. Kaneta's
visit in Perugia in October 2015.}  \\
Dipartimento di Matematica e Informatica\\
Universit\`a degli Studi di Perugia\\
Perugia (Italy)\\
\{gino,fernanda\}@dmi.unipg.it\\
Hitoshi Kaneta \\
Kyo-machi 77, Tsuyama, Okayama, Japan\\
hkaneta@marble.ocn.ne.jp\\
%Palace Mozu 301, Mozu-Ume-machi 3-34-8, \\ Kita-ku,  Sakai, 591-8032, Japan \\
%hkaneta@river.sannet.ne.jp \\ \\
}
\date{}

\maketitle

%\institute{ F. Pambianco ($\boxtimes$)\at Department of
%Mathematics and Informatics, Perugia University, Perugia,
%06123, Italy\\\email{fernanda@dmi.unipg.it}\\phone
%+39(075)5855006 \\fax~~~~~ +39(075)5855024\\
%\\
%H. Kaneta ($\boxtimes$)\at
%Palace Mozu 301, Mozu-Ume-machi 3-34-8, Kita-ku, Sakai, 591-8032, Japan\\
%\email{hkaneta@river.sannet.ne.jp}}
%\date{Received:  / Accepted: }

\begin{abstract} All  $\SSS{5}$-invariant nonsingular quartic surfaces are obtained. There exist no $\AAA{6}$-invariant
nonsingular quartic surfaces.
%to the Fermat curve of equation $x^d+y^d+z^d=0$. For some values
%of $d\leq 20$ this result has been obtained in \cite{enr,kmp1,kmp2,sem}.
%For the exceptional cases, $d=4,6$, the Klein quartic
%\cite{har} and the Wiman sextic \cite{doi} are respectively the
%uniquely determined maximally symmetric
%curves.\keywords{Algebraic curves  \and
%Automorphism groups \and Symmetric curves \and Cyclic subgroups}
% \subclass{14H45  \and 14N15}
\end{abstract}


\maketitle
\setcounter{section}{-1}
%%% section 0
\section{Introduction}
Let $V(f)$ be a nonsingular algebraic surface defined by a homogeneous polynomial $f(x,y,z,t)$ over an algebraically closed field $k$
of degree $d\geq 3$. The characteristic of the field is assumed to be zero.
The projective automorphism group $\Paut(V(f))$ coincides with the automorphism group $\Aut(V(f))$ unless $d=4$ \cite{mat}. Moreover, there exists a constant
$B_d>0$ such that $|\Paut(V(f))|\leq B_d$ for any nonsingular homogeneous polynomial $f$ of degree $d$.
Let $G$ be a group. If $\Paut(V(f))$ contains a subgroup isomorphic to $G$, $V(f)$ is said to be $G$-invariant.
When $G$ is a subgroup of $PGL_4(k)$, then $V(f)$ is said to be $G$-invariant if $\Paut(V(f))$ contains a subgroup conjugate to
$G$. The most symmetric nonsingular cubic surface is projectively equivalent to $V(x^3+y^3+z^3+t^3)$. Meanwhile the secondly most symmetric
nonsingular cubic surface $V(f)$ can be characterized in two ways: either $f$ is projectively equivalent to $x^2y+y^2z+z^2t+t^2x$ or $V(f)$ is
$\SSS{5}$-invariant \cite{dol}. Indeed, $\Aut(V(x^2y+y^2z+z^2t+t^2x))$ is isomorphic to $\SSS{5}$. Burnside conjectured that the most symmetric nonsingular
quartic surface is projectively equivalent to  $V(h)$, where $h=x^4+y^4+z^4+t^4+12xyzt$ so that $|\Paut(V(h))|=1920$ \cite[\S272]{bur}.

  In this paper $V(h)$ is shown to be $\SSS{5}$-invariant, and all $\SSS{5}$-invariant nonsingular quartic surfaces are given, up to projective equivalence.
%%It remains, however, to classify these surfaces and to determine the projective automorphism groups of these surfaces.
All $\AAA{5}$-invariant quartic surfaces, hence $\SSS{5}$-invariant  quartic surfaces as well, are given by Dolgachev in an intrinsic way \cite{dol2}. Unlike
\cite{dol2} we start with classifications of the faithful representations of $\AAA{5}$ and $\SSS{5}$ in $PGL_4(k)$ to get invariant nonsingular quartic
forms in four variables.
Besides, it is also shown that there
exist no $\AAA{6}$-invariant nonsingular quartic surfaces.

%%% section 1
\section{Preliminaries}

Let $\omega\in k^*$ be of order three, $\omega=(-1+i\sqrt{3})/2$. $M_{m,n}(k)$ stands for the set of all $m\times n$ matrices with entries in $k$.
By definition $M_n(k)=M_{n,n}(k)$, $GL_n(k)=\{A=[a_{ij}]\in M_n(k)\ :\ \det A\not=0\}$, and $PGL_n(k)=GL_n(k)/( E_n)$,
where $(E_n)$ is the subgroup $\{\lambda E_n\ : \ \lambda\in k^*\}$ ($E_n$ is the unit matrix in $GL_n(k)$). The $i$-th column vector of $E_4$ will be
denoted by $e_i$ ($i\in [1,4]$).
The coset $A(E_n)$ containing an $A\in GL_n(k)$ will be  denoted $(A)$. We denote by $k[x]$ the $k$-algebra of polynomials in $x=[x_1,...,x_n]$ over $k$.
For an $A\in GL_n(k)$ and $f\in k[x]$ we define a polynomial $f_A\in k[x]$ to be $f_A(x)=f(\sum \alpha_{1j}x_j,...,\sum \alpha_{nj}x_j)$, where
$A^{-1}=[\alpha_{ij}]$. As  is well known, the map $T_A:k[x]\rightarrow k[x]$ assigning $f_A$ to $f$ is a $k$-algebra isomorphism of $k[x]$ such that
$T_AT_B=T_{AB}$, that is, $(f_B)_A=f_{AB}$. In particular, if $S\in GL_n(k)$ and $B=S^{-1}AS$, then $f_A\sim f$ if and only if $(f_{S^{-1}})_B\sim f_{S^{-1}}$.

A homogeneous polynomial $f$ of degree $d\geq 1$ defines a projective algebraic set $$V(f)=\{(a)\in P^{n-1}: \ f(a)=0\}$$ of an $(n-1)-$dimensional
projective space $P^{n-1}$ over $k$. $V(f)$ is called a hypersurface of degree $d$.
Let $a=[a_1,...,a_n]\in k^n$, $(a)\in V(f)$, and $A\in GL_n(k)$. Then $(a)$ is a singular point of $V(f)$ if $f_{x_i}(a)=0$ for all $i$. If
$(a)$ is a nonsingular point of $V(f)$,  $V(\sum_{i=1}^n\gamma_ix_i)$ is the tangent plane to $V(f)$ at $(a)$, where $\gamma_i=f_{x_i}(a)$. Clearly
$(A):V(f)\rightarrow V(f_{A})$ is a bijection, and if $b=Aa$ with $(a)\in V(f)$, then $(f_A)_{x_j}(b)=\sum_{i=1}^n \gamma_i\alpha_{ij}$, where
$A^{-1}=[\alpha_{ij}]$. Consequently $(b)=(A)(a)$ is a nonsingular point of $V(f_A)$ if and only if $(a)$ is a nonsingular point of $V(f)$, and the tangent
plane of $V(f_A)$ at $(b)$ coincides with $(A)V(\sum_{i=1}^n \gamma_ix_i)$. In particular if $(a)$ is a nonsingular point of $V(f)$, $f_{A}\sim f$ and
$(A)(a)=(a)$, then $[f_{x_1}(a),...,f_{x_n}(a)]A\sim [f_{x_1}(a),...,f_{x_n}(a)]$. As is well know, $V(f)$ is irreducible if and only if $f=h^m$ for
some irreducible and homogeneous polynomial $h$ and some positive integer $m$. So we may assume that $f$ is irreducible if $V(f)$ is nonsingular.
Let $\Aut(f)=\{A\in GL_n(k):\ f_A=f\}$, $\Paut(f)=\{(A)\in PGL_n(k):\ f_A\sim f\}$, and $\Paut(V(f))=\{(A)\in PGL_n(k):\ (A)V(f)=V(f)\}$.
If $f$ is irreducible, then
$\Paut(V(f))=\Paut(f)=\Aut(f)/\langle (\varepsilon E_n)\rangle$, where $\ord(\varepsilon)=d$.

Let $G$ and $H$ be groups. A group homomorphism $\varphi$ of $H$ into $G$ is called a representation of $H$ in $G$. Two representations $\varphi$ and $\psi$
of $H$ in $G$ are said to be equivalent if there exists a $g\in G$ such that $\psi(h)=g^{-1}\varphi(h)g$ for any $h\in H$. A representation $\varphi$ is called
faithful if $\varphi$ is injective. We denote the symmetric group and alternating group of $n$ elements by $\SSS{n}$ and $\AAA{n}$ respectively.
Assume $n\geq 3$, and let $t_i=(i\ i+1)$ ($i\in [1,n-1]$).  Here $(ij)$ stands for a transposition.
Then $t_i^2=1$, $(t_it_{i+1})^3=1$, and $(t_it_j)^2=1$($|i-j|\geq 2$). Let $s_1=(123)$ and $s_{j}=(12)(j+1\ j+2)$ ($j\in [2,n-2]$). Then
$s_1^3=1$, $s_i^2=1$($i\in [2,n-2]$), $(s_{i-1}s_i)^3=1$($i\in [2,n-2]$), and $(s_is_j)^2=1$($|i-j|\geq 2$).
The following theorem \cite[chap.3,\S2]{suz} is due to Moore \cite{moo}.
%%%theorem 2.1
\begin{theorem}
 Let $G$ be a group.\\
$(1)$ There exists a faithful representation $\varphi$ of $\SSS{n}$ in $G$ such that $\tau_i=\varphi(t_i)$ $(i\in [1,n-1])$,
if and only if $\tau_i$ $(i\in [1,n-1])$ satisfy
\[
 \ord(\tau_i)=2,\ \ord(\tau_i\tau_{i+1})=3\ (i\in [1,n-2]),\ {\rm and}\ \ord(\tau_i\tau_j)=2\ (|i-j|\geq 2).
\]
$(2)$ There exists a faithful representation $\varphi$ of $\AAA{n}$ in $G$ such that $\sigma_i=\varphi(s_i)$ $(i\in [1,n-2])$,
if and only if $\sigma_i$ $(i\in [1,n-2])$ satisfy $\ord(\sigma_1)=3$,
\[
 \ord(\sigma_i)=2\ (i\in [2,n-2]),\ \ord(\sigma_{i-1}\sigma_i)=3\ (i\in [2,n-2]),\ {\rm and}\ \ord(\sigma_i\sigma_j)=2\ (|i-j|\geq 2).
\]
\end{theorem}

%In $\S 1$ it will be shown that the Fermat surface $S=V(x^3+y^3+z^3+t^3)$ is the unique cubic nonsingular surface, up to projective equivalence, such that
%$\Aut(S)$ is isomorphic to  $(\Z_3)^3\times_s \SSS{4}$.
%In $\S 2$ it will be shown that a surface $S'=V(x^2t+y^2z+z^2t+t^2x)$ is the unique
%cubic nonsingular surface such that $\Aut(S')$ is isomorphic to $\SSS{5}$.
%In $\S 3$ we classify the projective automorphism groups $\Paut(C_3)$ of nonsingular cubic curves.

%%%%% section 2
\section{Faithful representations of $\AAA{5}$, $\SSS{5}$ and $\AAA{6}$ in $PGL_4(k)$}
All faithful representations of $\AAA{5}$, $\SSS{5}$ and $\AAA{6}$ in $PGL_4(k)$ are found by Maschke \cite{mas} up to equivalence.
We describe these representations for our later use. We begin with the faithful representations $\varphi_i$($i\in [1,5]$) of $\AAA{5}$.
Let $s_1=(123)$, $s_2=(12)(34)$ and $s_3=(12)(45)$, and $(Q_{ij})=\varphi_i(s_j)$($j\in [1,3]$), where $Q_{ij}\in GL_4(k)$ are given as follows.
\begin{eqnarray*}
Q_{11}&=&\left[\begin{array}{cccc}
            1&0&0&0\\
            0&1&0&0\\
            0&0&\omega&0\\
            0&0&0&\omega^2\end{array}\right],\
Q_{12}=\left[\begin{array}{cccc}
            1&0&0&0\\
            0&-\frac{1}{3}&\frac{2}{3}&\frac{2}{3}\\
            0&\frac{2}{3}&-\frac{1}{3}&\frac{2}{3}\\
            0&\frac{2}{3}&\frac{2}{3}&-\frac{1}{3}\end{array}\right],\
Q_{13}=\left[\begin{array}{cccc}
            1&0&0&0\\
            0&-1&0&0\\
            0&0&0&\lambda_+\\
            0&0&\lambda_-&0\end{array}\right],\\
Q_{21}&=&Q_{11},\ \ \
Q_{22}=Q_{12},\ \ \
Q_{23}=\left[\begin{array}{cccc}
            -\frac{1}{4}&\frac{\sqrt{15}}{4}&0&0\\
            \frac{\sqrt{15}}{4}&\frac{1}{4}&0&0\\
            0&0&0&1\\
            0&0&1&0\end{array}\right],\\
Q_{31}&=&Q_{11},\ \ \
Q_{32}=\left[\begin{array}{cccc}
            \frac{1}{\sqrt{3}}&0&0&\frac{\sqrt{2}}{\sqrt{3}}\\
            0&-\frac{1}{\sqrt{3}}&\frac{\sqrt{2}}{\sqrt{3}}&0\\
            0&\frac{\sqrt{2}}{\sqrt{3}}&\frac{1}{\sqrt{3}}&0\\
            \frac{\sqrt{2}}{\sqrt{3}}&0&0&-\frac{1}{\sqrt{3}}\end{array}\right],\
%\frac{1}{\sqrt{3}}\left[\begin{array}{cccc}
%            1&0&0&\sqrt{2}\\
%            0&-1&\sqrt{2}&0\\
%            0&\sqrt{2}&1&0\\
%            \sqrt{2}&0&0&-1\end{array}\right],\
Q_{33}=\left[\begin{array}{cccc}
            \frac{\sqrt{3}}{2}&\frac{1}{2}&0&0\\
            \frac{1}{2}&-\frac{\sqrt{3}}{2}&0&0\\
            0&0&0&1\\
            0&0&1&0\end{array}\right],\\
Q_{41}&=&\left[\begin{array}{cccc}
            \omega&0&0&0\\
            0&\omega&0&0\\
            0&0&\omega^2&0\\
            0&0&0&\omega^2\end{array}\right],\
Q_{42}=\left[\begin{array}{cccc}
            \frac{1}{\sqrt{3}}&0&0&\frac{\sqrt{2}}{\sqrt{3}}\\
            0&\frac{1}{\sqrt{3}}&\frac{\sqrt{2}}{\sqrt{3}}&0\\
            0&\frac{\sqrt{2}}{\sqrt{3}}&-\frac{1}{\sqrt{3}}&0\\
            \frac{\sqrt{2}}{\sqrt{3}}&0&0&-\frac{1}{\sqrt{3}}\end{array}\right],\
%\frac{1}{\sqrt{3}}\left[\begin{array}{cccc}
%            1&0&0&\sqrt{2}\\
%            0&1&\sqrt{2}&0\\
%            0&\sqrt{2}&-1&0\\
%            \sqrt{2}&0&0&-1\end{array}\right],\
Q_{43}=\left[\begin{array}{cccc}
            0&0&0&\nu_{+}\\
            0&0&\nu_{+}&0\\
            0&\nu_{-}&0&0\\
            \nu_{-}&0&0&0\end{array}\right],\\
Q_{51}&=&Q_{41},\
Q_{52}=Q_{42},\
Q_{53}=\left[\begin{array}{cccc}
            0&0&0&\nu_{+}\\
            0&0&\nu_{-}&0\\
            0&\nu_{+}&0&0\\
            \nu_{-}&0&0&0\end{array}\right],
\end{eqnarray*}
where $\omega=(-1+i\sqrt{3})/2$, $\lambda_\pm=(-1 \pm i\sqrt{15})/4$, $\sqrt{15}=\sqrt{3}\sqrt{5}$, and $\nu_{\pm}=(\sqrt{3}\pm i\sqrt{5})/2\sqrt{2}$.

The representations $\varphi_i$ and $\varphi_j$ ($i\not=j$) are not equivalent, and any faithful representation of $\AAA{5}$ in $PGL_4(k)$ is
equivalent to one of $\varphi_i$ \cite{mas}. If we write
\begin{eqnarray*}
&&\varphi_i=\varphi_{i,\sqrt{3},\sqrt{5}}\ \ (i\in [1,2]),\\
&&\varphi_3=\varphi_{3,\sqrt{2},\sqrt{3}},\\
&&\varphi_i=\varphi_{i,\sqrt{2},\sqrt{3},\sqrt{5}}\ \ (i\in [4,5]),
\end{eqnarray*}
then $\varphi_{i,\pm\sqrt{3},\pm\sqrt{5}}$, $\varphi_{3,\sqrt{2},\sqrt{3}}$ and $\varphi_{i,\pm\sqrt{2},\pm\sqrt{3},\pm\sqrt{5}}$ also are
faithful representations of $\AAA{5}$.
%% lemma 2.1
\begin{lemma}
$(1)$ The representation $\varphi_{1,\pm\sqrt{3},\pm\sqrt{5}}$ is equivalent to one of the two representations $\varphi_{1,\sqrt{3},\sqrt{5}}$ and
$\varphi_{1,\sqrt{3},-\sqrt{5}}$, which are not equivalent. \\
$(2)$ The representation $\varphi_{2,\pm\sqrt{3},\pm\sqrt{5}}$ is equivalent to $\varphi_{2,\sqrt{3},\sqrt{5}}$.\\
$(3)$ The representation $\varphi_{3,\pm\sqrt{2},\pm\sqrt{3}}$ is equivalent to $\varphi_{3,\sqrt{2},\sqrt{3}}$.\\
$(4)$ The representation $\varphi_{4,\sqrt{2},\sqrt{3},\sqrt{5}}$ is equivqlent to one of the two representations
$\varphi_{4,\sqrt{2},\sqrt{3},\sqrt{5}}$ and $\varphi_{4,\sqrt{2},\sqrt{3},-\sqrt{5}}$.\\
$(5)$ The representation $\varphi_{5,\pm\sqrt{2},\pm\sqrt{3},\pm\sqrt{5}}$ is equivalent to $\varphi_{5,\sqrt{2},\sqrt{3},\sqrt{5}}$.
\end{lemma}
\Proof
(1) Let $\psi_1=\psi_{1,\sqrt{3},\sqrt{5}}$ be a representation of $\AAA{5}$ in $GL_4(k)$ such that
$\psi(s_j)=Q_{1j}$ ($j\in [1,3]$). Note that $\varphi_1=\pi\circ\psi_1$,
where $\pi$ is the canonical projection of $GL_4(k)$ onto $PGL_k(4)$. For $T=[e_1,e_2,e_4,e_3]$ we have
$T^{-1}\psi_{1}(s)T=\psi_{1,-\sqrt{3},\sqrt{5}}(s)$, for any $s=s_j$, hence for any $s\in \AAA{5}$. Therefore, $\varphi_{1,\sqrt{3},\sqrt{5}}$
 (resp. $\varphi_{1,\sqrt{3},-\sqrt{5}}$) is equivalent to $\varphi_{1,-\sqrt{3},\sqrt{5}}$ (resp. $\varphi_{1,-\sqrt{3},-\sqrt{5}}$). It is
easy to see that there exists no $U\in GL_4(k)$ such that $U^{-1}Q_{1j}(\sqrt{3})U\sim Q_{1j}(-\sqrt{3})$ ($j\in [1,3]$).\\
(2) Let $\psi_2=\psi_{2,\sqrt{3},\sqrt{5}}$ be a representation of $\AAA{5}$ in $GL_4(k)$ such that
$\psi(s_j)=Q_{2j}$ ($j\in [1,3]$). Note that $\varphi_2=\pi\circ\psi_2$. For $T=[-e_1,e_2,e_3,e_4]$ and $U=[-e_1,e_2,e_4,e_3]$ we have
\begin{eqnarray*}
&&T^{-1}Q_{21}(\sqrt{3})T=Q_{21}(\sqrt{3}),\ T^{-1}Q_{22}T=Q_{22},\ T^{-1}Q_{23}(\sqrt{3},\sqrt{5})T=Q_{23}(\sqrt{3},-\sqrt{5}),\\
&&U^{-1}Q_{21}(\sqrt{3})U=Q_{21}(-\sqrt{3}),\ U^{-1}Q_{22}U=Q_{22},\ U^{-1}Q_{23}(\sqrt{3},\sqrt{5})U=Q_{23}(-\sqrt{3},\sqrt{5}).
\end{eqnarray*}
(3) For $T=[-e_1,-e_2,e_3,e_4]$ we have $T^{-1}Q_{32}(\sqrt{2},\sqrt{3})T=Q_{32}(-\sqrt{2},\sqrt{3})$, hence $\varphi_{3,\sqrt{2},\sqrt{3}}$
and $\varphi_{3,-\sqrt{2},\sqrt{3}}$ are equivalent. For $U=[e_2,e_1,e_4,e_3]$ we have
\begin{eqnarray*}
&&U^{-1}Q_{31}(\sqrt{3})U=Q_{31}(-\sqrt{3}),\ U^{-1}Q_{32}(\sqrt{2},\sqrt{3})U=Q_{32}(-\sqrt{2},-\sqrt{3}),\\
&&U^{-1}Q_{33}(\sqrt{3})U=Q_{33}(-\sqrt{3}),
\end{eqnarray*}
so that $\varphi_{3,\sqrt{2},\sqrt{3}}$ and $\varphi_{3,-\sqrt{2},-\sqrt{3}}$ are equivalent.\\
(4) For $T=[e_4,e_3,e_2,e_1]$ we have
\begin{eqnarray*}
&&T^{-1}Q_{4,1}(\sqrt{3})T=Q_{41}(-\sqrt{3}),\ T^{-1}Q_{42}(\sqrt{2},\sqrt{3})T=Q_{42}(-\sqrt{2},-\sqrt{3}),\\
&&T^{-1}Q_{43}(\sqrt{2},\sqrt{3},\sqrt{5})T=-Q_{43}(-\sqrt{2},-\sqrt{3},\sqrt{5}),
\end{eqnarray*}
hence $\varphi_{4,\varepsilon_1\sqrt{2},\varepsilon_2\sqrt{3},\varepsilon_3\sqrt{5}}$ ($\varepsilon_i=\{\pm\}$) and
$\varphi_{4,-\varepsilon_1\sqrt{2},-\varepsilon_2\sqrt{3},\varepsilon_3\sqrt{5}}$ are equivalent.
For $T=[e_4,e_3,-e_2,-e_1]$ we have
\begin{eqnarray*}
&&T^{-1}Q_{41}(\sqrt{3})T=-Q_{41}(-\sqrt{3}),\ T^{-1}Q_{42}(\sqrt{2},\sqrt{3})T=-Q_{42}(\sqrt{2},\sqrt{3}),\\
&&T^{-1}Q_{43}(\sqrt{3},\sqrt{5})T=-Q_{43}(-\sqrt{3},\sqrt{5}),
\end{eqnarray*}
hence $\varphi_{4,\varepsilon_1\sqrt{2},\varepsilon_2\sqrt{3},\varepsilon_3\sqrt{5}}$ and
$\varphi_{4,\varepsilon_1\sqrt{2},-\varepsilon_2\sqrt{3},\varepsilon_3\sqrt{5}}$ are equivalent. However,
$\varphi_{4,\sqrt{2},\sqrt{3},\sqrt{5}}$ and $\varphi_{4,\sqrt{2},\sqrt{3},-\sqrt{5}}$ are not equivalent, for any $U\in GL_4(k)$ such that
\begin{eqnarray*}
&&U^{-1}Q_{41}(\sqrt{3})U\sim Q_{41}(\sqrt{3}),\ U^{-1}Q_{43}(\sqrt{2},\sqrt{3},\sqrt{5})U\sim Q_{43}(\sqrt{2},\sqrt{3},-\sqrt{5})
\end{eqnarray*}
does not satisfy $U^{-1}Q_{42}(\sqrt{2},\sqrt{3})U\sim Q_{42}(\sqrt{2},\sqrt{3})$.\\
(5) For $T=[e_3,e_4,e_1,e_2]$, $U=[e_2,e_1,e_4,e_3]$ and $V=[e_4,e_3,-e_2,-e_1]$  we have
\begin{eqnarray*}
&&T^{-1}Q_{51}(\sqrt{3})T=Q_{51}(-\sqrt{3}),\ T^{-1}Q_{52}(\sqrt{2},\sqrt{3})T=Q_{52}(-\sqrt{2},-\sqrt{3}),\\
&&T^{-1}Q_{53}(\sqrt{2},\sqrt{3},\sqrt{5})T=Q_{53}(-\sqrt{2},-\sqrt{3},-\sqrt{5}),\\
&&U^{-1}Q_{51}(\sqrt{3})U=Q_{51}(\sqrt{3}),\ U^{-1}Q_{52}(\sqrt{2},\sqrt{3})U=Q_{52}(\sqrt{2},\sqrt{3}),\\
&&U^{-1}Q_{53}(\sqrt{2},\sqrt{3},\sqrt{5})U=Q_{53}(\sqrt{2},\sqrt{3},-\sqrt{5}),\\
&&V^{-1}Q_{51}(\sqrt{3})V=Q_{51}(-\sqrt{3}),\ V^{-1}Q_{52}(\sqrt{2},\sqrt{3})V=Q_{52}(\sqrt{2},-\sqrt{3}),\\
&&V^{-1}Q_{53}(\sqrt{2},\sqrt{3},\sqrt{5})V=Q_{53}(\sqrt{2},-\sqrt{3}, \sqrt{5}).
\end{eqnarray*}
Consequently any representation $\varphi_{5,\pm\sqrt{2},\pm\sqrt{3},\pm\sqrt{5}}$ is equivalent to $\varphi_{5,\sqrt{2},\sqrt{3},\sqrt{5}}$.\\

 We proceed to describe the faithful representations $\Phi_i$ ($i\in [1,3]$) of $\SSS{5}$ in $PGL_4(k)$. Let $t_j=(j\ j+1)\in \SSS{5}$ ($j\in [1,4]$).
Recall $\omega=(-1+i\sqrt{3})/2$. Let $\Phi_i(s_j)=(R_{ij})$ ($j\in [1,3]$), and $\Phi_i(t_1)=(R_{i4})$, where $R_{ij}\in GL_4(k)$ are given as follows.
\begin{eqnarray*}
R_{11}&=&\left[\begin{array}{cccc}
            1&0&0&0\\
            0&1&0&0\\
            0&0&\omega&0\\
            0&0&0&\omega^2\end{array}\right],\
R_{12}=\left[\begin{array}{cccc}
            1&0&0&0\\
            0&-\frac{1}{3}&\frac{2}{3}&\frac{2}{3}\\
            0&\frac{2}{3}&-\frac{1}{3}&\frac{2}{3}\\
            0&\frac{2}{3}&\frac{2}{3}&-\frac{1}{3}\end{array}\right],\
R_{13}=\left[\begin{array}{cccc}
            -\frac{1}{4}&\frac{\sqrt{15}}{4}&0&0\\
            \frac{\sqrt{15}}{4}&\frac{1}{4}&0&0\\
            0&0&0&1\\
            0&0&1&0\end{array}\right],\\
R_{14}&=&\left[\begin{array}{cccc}
                             1&0&0&0\\
                             0&1&0&0\\
                             0&0&0&1\\
                             0&0&1&0\end{array}\right].
\end{eqnarray*}
\begin{eqnarray*}
R_{21}&=&\left[\begin{array}{cccc}
            1&0&0&0\\
            0&1&0&0\\
            0&0&\omega&0\\
            0&0&0&\omega^2\end{array}\right],\
R_{22}=\frac{1}{\sqrt{3}}\left[\begin{array}{cccc}
            1&0&0&\sqrt{2}\\
            0&-1&\sqrt{2}&0\\
            0&\sqrt{2}&1&0\\
            \sqrt{2}&0&0&-1\end{array}\right],\
R_{23}=\left[\begin{array}{cccc}
            \frac{\sqrt{3}}{2}&\frac{1}{2}&0&0\\
            \frac{1}{2}&-\frac{\sqrt{3}}{2}&0&0\\
            0&0&0&1\\
            0&0&1&0\end{array}\right],\\
R_{24}&=&\left[\begin{array}{cccc}
                             0&1&0&0\\
                             -1&0&0&0\\
                             0&0&0&1\\
                             0&0&-1&0\end{array}\right].
\end{eqnarray*}
\begin{eqnarray*}
R_{31}&=&\left[\begin{array}{cccc}
            \omega&0&0&0\\
            0&\omega&0&0\\
            0&0&\omega^2&0\\
            0&0&0&\omega^2\end{array}\right],\
R_{32}=\frac{1}{\sqrt{3}}\left[\begin{array}{cccc}
            1&0&0&\sqrt{2}\\
            0&1&\sqrt{2}&0\\
            0&\sqrt{2}&-1&0\\
            \sqrt{2}&0&0&-1\end{array}\right],\
R_{33}=\left[\begin{array}{cccc}
            0&0&0&\nu_+\\

            0&0&\nu_-&0\\

            0&\nu_+&0&0\\

            \nu_-&0&0&0\end{array}\right],\\
R_{34}&=&\left[\begin{array}{cccc}

                             0&0&1&0\\
                             0&0&0&-1\\
                             1&0&0&0\\
                             0&-1&0&0\end{array}\right],
\end{eqnarray*}
where $\nu_{\pm}=(\sqrt{3}\pm i\sqrt{5})/2\sqrt{2}$.  Note that $R_{34}$ is wrongly defined to be $[-e_4,e_3,-e_2,e_1]$ in \cite[p.296]{mas}.
The representations $\Phi_i$ and $\Phi_j$ ($i\not=j$) are not equivalent, and any faithful representation of $\SSS{5}$ in $PGL_4(k)$ is equivalent to one of
$\Phi_i$ ($i\in [1,3]$)\cite{mas}. To be more precise, if we write $\Phi_1=\Phi_{1,\sqrt{3},\sqrt{5}}$, $\Phi_2=\Phi_{2,\sqrt{2},\sqrt{3}}$, and
$\Phi_3=\Phi_{3,\sqrt{2},\sqrt{3},\sqrt{5}}$, then $\Phi_{1,\pm\sqrt{3},\pm\sqrt{5}}$, $\Phi_{2,\pm\sqrt{2},\pm\sqrt{3}}$ and
$\Phi_{3,\pm\sqrt{2},\pm\sqrt{3},\pm\sqrt{5}}$ also are faithful representations of $\SSS{5}$.
%% lemmma 2.2
\begin{lemma}
$(1)$ The representations  $\Phi_{1,\pm\sqrt{3},\pm\sqrt{5}}$ are equivalent to $\Phi_{1,\sqrt{3},\sqrt{5}}$.\\
$(2)$ The representations $\Phi_{2,\pm\sqrt{2},\pm\sqrt{3}}$ are equivalent to $\Phi_{2,\sqrt{2},\sqrt{3}}$.\\
$(3)$ The representations $\Phi_{3,\pm\sqrt{2},\pm\sqrt{3},\pm\sqrt{5}}$ are equivalent to $\Phi_{3,\sqrt{2},\sqrt{3},\sqrt{5}}$.
\end{lemma}
\Proof
Note that $\Phi_1$, $\Phi_2$ and $\Phi_3$ are extensions of $\varphi_2$, $\varphi_3$ and $\varphi_5$, respectively.
(1) Let $\Psi_1$ be the representation of $\SSS{5}$ in $GL_4(k)$ such that $\Psi_1(s_i)=R_{1i}$ ($i\in [1,3]$) and  $\Psi_1(t_1)=R_{14}$.
Then we can easily verify that $\Psi_1=\Psi_{1,\sqrt{3},\sqrt{5}}$ is equivalent to the irreducible representation denoted by $V$ on \cite[p.28]{ful}. Since
$X^{-1}R_{14}X=R_{14}$ for any $X\in \{T,U\}$ defined in the proof of Lemmma 2.1(2), the representations $\Psi_{1,\pm\sqrt{3},\pm\sqrt{5}}$
are equivalent. Thus (1) follows, for $\Phi_1=\pi\circ\Psi_1$.
(2) holds, for $X^{-1}R_{24}X\sim R_{24}$ for any $X\in \{T,U\}$ defined in the proof of Lemma 2.1(3). Simillarly (3) holds, for
$X^{-1}R_{34}X\sim R_{34}$ for any $X\in \{T,U,V\}$ defined in the proof of Lemma 2.1(5).\\

Finally we describe all faithful representations of $\AAA{6}$ in $PGL_4(k)$. Define matrices $Q_{ij}\in GL_4(k)$ ($i\in [6,7],j\in [1,4]$) as follows;
$Q_{6j}=Q_{3j}$, $Q_{7j}=Q_{5j}$ ($j\in [1,3]$) and
\begin{eqnarray*}
&&Q_{64}=\left[\begin{array}{cccc}
0&1&0&0\\
1&0&0&0\\
0&0&0&-1\\
0&0&-1&0\end{array}\right],\
Q_{74}=\frac{1}{\sqrt{5}}\left[\begin{array}{cccc}
0&0&\sqrt{3}&i\sqrt{2}\\
0&0&-i\sqrt{2}&-\sqrt{3}\\
\sqrt{3}&i\sqrt{2}&0&0\\
-i\sqrt{2}&-\sqrt{3}&0&0\end{array}\right].
\end{eqnarray*}
By Theorem 1.1 there exist faithful representations $\varphi_i$ ($i\in [6,7]$) of $\AAA{6}$ in $PGL_4(k)$ such that $\varphi_i(s_j)=(Q_{ij})$ ($j\in [1,4]$).
They are not equivalent, and any faithful representation of $\AAA{6}$ is equivalent to one of them \cite{mas}. To be more precise, if we denote
the representations $\varphi_6$ and $\varphi_7$ by $\varphi_{6,\sqrt{2},\sqrt{3}}$ and $\varphi_{7,\sqrt{2},\sqrt{3},\sqrt{5}}$, then
$\varphi_{6,\pm\sqrt{2},\pm\sqrt{3}}$ and $\varphi_{7,\pm\sqrt{2},\pm\sqrt{3},\pm\sqrt{5}}$ also are representations.

%% lemma 2.3
\begin{lemma}
$(1)$ The representations $\varphi_{6,\pm\sqrt{2},\pm\sqrt{3}}$ are equivalent to $\varphi_{6,\sqrt{2},\sqrt{3}}$.\\
$(2)$ The representations $\varphi_{7,\pm\sqrt{2},\pm\sqrt{3},\pm\sqrt{5}}$ are equivalent to $\varphi_{7,\sqrt{2},\sqrt{3},\sqrt{5}}$.
\end{lemma}
\Proof
(1) Let $T$ and $U$ be as in the proof of Lemma 2.1(3). Then $X^{-1}Q_{64}X\sim Q_{64}$ for any $X\in \{T,U\}$. (2) Let $T$, $U$, and $V$ be as
in the proof of Lemma 2.1(5). Then $X^{-1}Q_{74}X\sim Q_{74}$ for any $X\in \{T,U,V\}$.\\
%%%%% section 3
\section{$\AAA{5}$-invariant quartic surfaces}
As we have seen in the previous section, all faithful representations of  $\AAA{5}$ in $PGL_k(4)$ are
 $\varphi_{1,\sqrt{3},\pm\sqrt{5}}$, $\varphi_{2,\sqrt{3},\sqrt{5}}$, $\varphi_{3,\sqrt{2},\sqrt{3}}$, $\varphi_{4,\sqrt{2},\sqrt{3},\pm\sqrt{5}}$
and $\varphi_{5,\sqrt{2},\sqrt{3},\sqrt{5}}$ up to equivalence. We will denote these representations by
$\varphi_1$, $\varphi_2$, $\varphi_3$, $\varphi_4$, $\varphi_5$, respectively. Let $\AAA{5}(i)=\varphi_{i}(\AAA{5})$ ($i\in [1,5]$).

\begin{lemma} \label{lem 3.1} No nonsingular quartic surface is $\AAA{5}(1)$-invariant.
\end{lemma}
\Proof
Suppose a quartic form $f$ to be $\AAA{5}(1)$-invariant. In particular $f_{Q_{11}^{-1}}=\omega^i f$. According as
$i=0$, $i=1$ or $i=2$, $f$ belongs to
\begin{eqnarray*}
&&\langle x^4,y^4,x^3y,y^3x,z^3x,z^3y,t^3x,t^3y,x^2y^2,z^2t^2,x^2zt,y^2zt,xyzt \rangle,\\
&&\langle z^4,x^3z,y^3z,t^3z,x^2t^2,y^2t^2,x^2yz,y^2xz,z^2xt,z^2yt,t^2xy \rangle,\\
&&\langle t^4,x^3t,y^3t,z^3t,x^2z^2,y^2z^2,x^2yt,y^2xt,z^2xy,t^2xz,t^2yz \rangle.
\end{eqnarray*}
  If $i=2$, then $f$ does not contain neither $t^4$ nor $x^3t$, for $f_{Q_{13}^{-1}}\sim f$, hence $V(f)$ is
singular at $(1,0,0,0)$. Similarly, $V(f)$ is singular at $(1,0,0,0)$ if $i=1$. Assume $i=0$ and
\[
 f=ax^4+bx^3y+x^2(c_1y^2+c_2zt)+ x(d_1y^3+d_2z^3+d_3t^3+d_4yzt)+e_1y^4+e_2y^2zt+e_3z^3y+e_4t^3y+e_5z^2t^2.
\]
   Since $f_{Q_{13}^{-1}}=\pm f$, we obtain $d_2\lambda_+^3=\pm d_3$ and $-e_3\lambda_+^3=\pm e_4$. Besides,
$f_{Q_{12}^{-1}}$ can contain none of monomials $xy^2z,\ xy^2t,\ z^4,\ t^4$, hence $d_2=d_3$ and $e_3=e_4$,
namely $d_2=d_3=e_3=e_4=0$. Consequently $V(f)$ is singular at $(0,0,0,1)$, for $f$ contains none of monomials
$t^3x,\ t^3y,\ t^3z,\ t^4$.\\

We denote the following quartic forms by $g_0$ and $g_1$, respectively.
\begin{eqnarray*}
&&-x^4+2\sqrt{15}(y^3+z^3+t^3)x+10(z^3+t^3)y+13x^2y^2+5z^2t^2+(26x^2+20y^2)zt-6\sqrt{15}xyzt,\\
&&x^4+y^4+2x^2y^2+4z^2t^2+4(x^2+y^2)zt.
\end{eqnarray*}

\begin{lemma} \label{lem 3.2}
An $\AAA{5}(2)$-invariant nonsingular quartic form takes the form $g_0+\lambda g_1$
where $\lambda\in k$.
\end{lemma}
\Proof
Clearly $V(g_1)$ is singular at $(0,0,1,0)$. We will show that $g_0$ is nonsingular. Suppose $V(g_0)$ is singular at $(x,y,z,t)$,
namely, $g_{0x}$, $g_{0y}$, $g_{0z}$, and $g_{0t}$ vanish there. If $t=z=0$, then $x=y=0$. If $t=0$ and $z\not=0$, then $g_{0z}=g_{0t}=0$
imply $x=y=0$, hence $z=0$ for $g_{0x}=0$, a contradiction. Thus $t\not=0$. To see $x\not=0$, assume $x=0$. Since $g_{0x}=0$, it follows that
$0=y^3+z^3+t^3-3yzt$. In addition $z^3-t^3=0$, for $zg_{0z}-tg_{0t}=30y(z^3-t^3)$ and $y\not=0$.
Now $t=z\omega^i\not=0$ ($i\in [0,2]$), and $y=-\omega^{-i}z/2$ for $g_{0y}=0$ so that $0=y^3+z^3+t^3-3yzt=27z^3/8\not=0$, contradiction. Thus
$x\not=0$. Now we may assume $t=1$. Note that $(6\sqrt{15}x+30y)(z^3-1)=zg_{0z}-g_{0t}$. First suppose $6\sqrt{15}x+30y=0$, i.e. $y=-\sqrt{15}x/5$.
Then $0=5g_{0x}-\sqrt{15}g_{0y}=-2424x^3\not=0$, a contradiction. Next suppose $z^3=1$. Then $g_{0y}-yg_{0z}$ is equal to
\begin{eqnarray*}
&&6\sqrt{15}xz(2y^2z^2-yz-1)-10(2y^3z^3+30y^2z^2-30yz-20)\\
&=&\{6\sqrt{15}u-10(v+2)\}(2v+1)(v-1),
\end{eqnarray*}
where $u=xz$ and $v=yz$. Consequently, 1) $v=1$, 2) $v=-1/2$, or 3) $v=(3\sqrt{15}u -10)/5$. If 1) is the case, $zg_{0t}=0$ and $g_{0x}=0$ give a contradiction,
namely $u^2=-30/13$ and $u^2=39/2$. If 2) is the case, $zg_{0t}=0$ and $g_{0x}=0$ give $u(26u+9\sqrt{15})=0$ and $16u^3-182u-27\sqrt{15}=0$, a contradiction,
for a pair of algebraic equations $X(26X+9\sqrt{15})=0$ and $16X^3-182X-27\sqrt{15}=0$ has no common roots. Even if 3) is the case,
$zg_{0t}=0$ and $g_{0x}=0$ give $40u^2-6\sqrt{15}u+15=0$, and $584u^3-318\sqrt{15}u^2+795u+25\sqrt{15}=0$, a contradiction. Therefore, $V(g_0)$ has no
singular points.

Suppose a quartic form $f$ to be $\AAA{5}(2)$-invariant. If $f_{Q_{21}^{-1}}=\omega f$, then $V(f)$ is singular at
$(1,0,0,0)$, for $f$ cannot contain $x^3z$ because $f_{Q_{22}^{-1}}\sim f$, hence $f$ contains none of $x^4,\ x^3y,\ x^3z,\ x^3t$.
Similarly, if $f_{Q_{21}^{-1}}=\omega^2 f$, then $V(f)$ is singular at $(1,0,0,0)$.
We may assume that $f_{Q_{21}^{-1}}= f$, hence $f$ has the form
\[
 f=ax^4+bx^3y+x^2(c_1y^2+c_2zt)+ x(d_1y^3+d_2z^3+d_3t^3+d_4yzt)+e_1y^4+e_2y^2zt+e_3z^3y+e_4t^3y+e_5z^2t^2.
\]
It holds that  $f_{Q_{22}^{-1}}=\pm f$, for $Q_{22}^2=E_4$. We see $b=0$, for $f_{Q_{22}^{-1}}$ contains none of monomials
$x^3z,\ xy^2z,\ xy^2t,\ xz^2y,\ xz^2t$. Consequently, if $a=0$, $f$ is singular at $(1,0,0,0)$, hence $a\not=0$ and $f_{Q_{22}^{-1}}= f$.  Now $f_{Q_{22}^{-1}}= f$ is equivalent
to
\begin{eqnarray*}
&&(c_1y^2+c_2zt)_{Q_{22}^{-1}}=c_1y^2+c_2zt,\ (d_1y^3+d_2z^3+d_3t^3+d_4yzt)_{Q_{22}^{-1}}=d_1y^3+d_2z^3+d_3t^3+d_4yzt,\\
&&(e_1y^4+e_2y^2zt+e_3z^3y+e_4t^3y+e_5z^2t^2)_{Q_{22}^{-1}}=e_1y^4+e_2y^2zt+e_3z^3y+e_4t^3y+e_5z^2t^2.
\end{eqnarray*}
Namely, $c_2=2c_1$, $d_1=d_2=d_3=-d_4/3$, $e_3=e_4=(-4e_1+e_2)/2$, and $e_5=(12e_1+e_2)/4$.
We denote $c_1$ and $d_1$ by $c$ and $d$, respectively. Then $f_{Q_{23}^{-1}}$ takes the following form:
\begin{eqnarray*}
&&\frac{1}{256}(a+15c-15\sqrt{15}d+225e_1)x^4+\frac{1}{256}(90a+166c+42\sqrt{15}d+90e_1)x^2y^2\\
&&+\frac{1}{256}(-60\sqrt{15}a+28\sqrt{15}c+44d+4\sqrt{15}e_1)xy^3+\frac{1}{256}(225a+15c+\sqrt{15d}+e_1)y^4\\
&&+\frac{1}{256}(-4\sqrt{15}a-28\sqrt{15}c+180d+60\sqrt{15}e_1)x^3y\\
&&+\frac{1}{16}(2c+3\sqrt{15}d+15e_2)x^2zt+\frac{1}{16}(-4\sqrt{15}-42d+2\sqrt{15}e_2)xyzt\\
&&+\frac{1}{16}(30c-3\sqrt{15}d+e_2)y^2zt\\
&&+\frac{1}{4}(-d-2\sqrt{15}e_1+\frac{1}{2}\sqrt{15}e_2)xz^3+\frac{1}{4}(-d-2\sqrt{15}e_1+\frac{1}{2}\sqrt{15}e_2)xt^3\\
&&+\frac{1}{4}(\sqrt{15}d-2\sqrt{15}e_1+\frac{1}{2}\sqrt{15}e_2)yz^3+\frac{1}{4}(\sqrt{15}d-2\sqrt{15}e_1+\frac{1}{2}\sqrt{15}e_2)yt^3\\
&&+(3e_1+\frac{1}{4}e_2)z^2t^2.
\end{eqnarray*}
Note that $f_{Q_{23}^{-1}}=\pm f$, for $Q_{23}^2=E_4$. Suppose $f_{Q_{23}^{-1}}-f=0$. The coefficients of $x^4$, $x^2y^2$, $xy^3$, $y^4$ and $x^3y$ yield
\[
 -a-7c+3\sqrt{15}d+15e_1=0,\ \ -30c+13\sqrt{15}d+60e_1=0.
\]
The coefficients of $x^2zt$, $xyzt$ and $y^2zt$ yield  $-10c+\sqrt{15}d+5e_2=0$. The coefficients of $xz^3$ and $yz^3$ yield $-5d+\sqrt{15}e_3=0$.
Hence, we can show that $f_{Q_{23}^{-1}}=f$ if and only if
\begin{eqnarray*}
&&a=(-\sqrt{15}d+30e_1)/30,\ c=(13\sqrt{15}d+60e_1)/30,\\
&&e_2=(2\sqrt{15}d+12e_1)/3,\ e_3=(\sqrt{15}d)/3,\ e_5=(\sqrt{15}d+24e_1)/6.
\end{eqnarray*}
Setting $d=2\sqrt{15}a_0$ and $e_1=a_1$, we obtain $f$ such that $f_{Q_{23}^{-1}}=f=a_0g_0+a_1g_1$.
Finally assuming  $f_{Q_{23}^{-1}}+f=0$ we will show that $f=0$.
The coefficients of $x^2zt$, $xyzt$ and $y^2zt$ yield  $2c+e_2=0$, and $-3\sqrt{15}d+2e_2=0$.
The coefficients of $xz^3$ and $yz^3$ yield $\sqrt{15}d+5e_3=0$.
%% Since $e_3=(-4e_1+e_2)/2$, we have $e_1=(19\sqrt{15}d)/2$.
In addition the coefficients of $x^4$, $x^2y^2$, $xy^3$, $y^4$ and $x^3y$ yield
\begin{eqnarray*}
&&257a+15c-15\sqrt{15}d+225e_1=0, \\
&&90a+422c+42\sqrt{15}d+90e_1=0,\\
&&-15a+7c+20\sqrt{15}d+e_1=0,\\
&&225a+15c+\sqrt{15}d+257e_1=0,\\
&&-a-7c+3\sqrt{15}d+15e_1=0.
\end{eqnarray*}
It is easy to see that the second, the third and the 5th equalities together with the equality $32a-16\sqrt{15}d-32e_1=0$,
which is the difference of the first and the 4th equalities, imply $a=c=d=e_1=0$, hence $e_2=e_3=e_4=e_5=0$. Now $f=0$ follows.\\

\begin{lemma} \label{lem 3.3}
 Any $\AAA{5}(3)$-invariant nonsingular quartic form $f$ satisfies
$f_{Q_{31}^{-1}}=f$, $f_{Q_{32}^{-1}}f=f$, $f_{Q_{33}^{-1}}=f$, and has the form
\begin{eqnarray*}
&&\frac{\sqrt{3}}{6}(b_1+2b_2)x^4-\frac{\sqrt{3}}{6}(2b_1+b_2)y^4+b_1x^3y+b_2xy^3\\
&&+\sqrt{2}b_2x z^3+\frac{\sqrt{2}}{\sqrt{3}}(2b_1+b_2)yz^3+\frac{\sqrt{2}}{\sqrt{3}}(b_1+2b_2)xt^3-\sqrt{2}b_1yt^3\\
&&+\frac{\sqrt{3}}{2}(-b_1+b_2)x^2y^2+\frac{\sqrt{3}}{2}(-b_1+b_2)z^2t^2\\
&&+3b_1x^2zt+3b_2y^2zt+\sqrt{3}(b_1-b_2)xyzt.
\end{eqnarray*}
\end{lemma}
\Proof
Let $f$ be an  $\AAA{5}(3)$-invariant nonsingular quartic form. As we noted in the proof of Lemma 3.2, $f_{Q_{31}^{-1}}= f$ and
$f$ has the form
\[
 a_1x^4+a_2y^4+b_1x^3y+b_2y^3x+b_3z^3x+b_4z^3y+b_5t^3x+b_6t^3y+c_1x^2y^2+c_2z^2t^2+d_1x^2zt+d_2y^2zt+exyzt.
\]
Since $Q_{32}^2=E_4$, we have $f_{Q_{32}^{-1}}=(-1)^jf$, where $j\in [0,1]$. Thus we obtain
\begin{eqnarray*}
&&(a_1+2\sqrt{2}b_5)/9=(-1)^ja_1,\ (8\sqrt{2}a_1+5b_5)/9=(-1)^jb_5,\\
&&(a_2-2\sqrt{2}b_4)/9=(-1)^ja_2,\ (-8\sqrt{2}a_2+5b_4)/9=(-1)^jb_4,\\
&&(-b_1-2\sqrt{2}b_6+2d_1)/9=(-1)^jb_1,\ (-2\sqrt{2}b_1+b_6-2\sqrt{2}d_1)/9=(-1)^jb_6,\\
&&(6b_1-6\sqrt{2}b_6+3d_1)/9=(-1)^jd_1,\\
&&(-b_2+2\sqrt{2}b_3+2d_2)/9=(-1)^jb_2,\ (2\sqrt{2}b_2+b_3+2\sqrt{2}d_2)/9=(-1)^jb_3,\\
&&(6b_2+6\sqrt{2}b_3+3d_2)/9=(-1)^jd_2,\\
&&(c_1+4c_2-2e)/9=(-1)^jc_1,\ (4c_1+c_2-2e)/9=(-1)^jc_2,\ (-8c_1-8c_2+e)/9=(-1)^je,
\end{eqnarray*}
together with
\begin{eqnarray*}
&&2\sqrt{2}a_1-b_5=0,\ 2\sqrt{2}a_2+b_4=0,\\
&&\sqrt{2}b_1+b_6=0,\ 3b_1-d_1=0,\ \sqrt{2}b_2-b_3=0,\ 3b_2-d_2=0,\ c_1=c_2,\ e=-2c_1.
\end{eqnarray*}
Note that if $i=1$, then $f=0$, and that $f_{Q_{32}^{-1}}=f$ holds  if and only if
\begin{eqnarray*}
&&2\sqrt{2}a_1-b_5=0,\ 2\sqrt{2}a_2+b_4=0,\\
&&\sqrt{2}b_1+b_6=0,\ 3b_1-d_1=0,\ \sqrt{2}b_2-b_3=0,\ 3b_2-d_2=0,\ c_1=c_2,\ e=-2c_1.
\end{eqnarray*}
Similarly $f_{Q_{33}^{-1}}=(-1)^jf$ ($j\in [0,1]$).  In order to describe this equality concretely we introduce matrices $W_2$, $W_3$ and $W_5$ as follows:
\begin{eqnarray*}
W_2&=&\frac{1}{2}\left[\begin{array}{cc}
            \sqrt{3}&1\\
            1&-\sqrt{3}
            \end{array}\right], \
W_3=\frac{1}{4}\left[\begin{array}{ccc}
            3&1&\sqrt{3}\\
            1&3&-\sqrt{3}\\
            2\sqrt{3}&-2\sqrt{3}&-2
            \end{array}\right],\\
W_5&=&\frac{1}{16}\left[\begin{array}{ccccc}
            9&1&3\sqrt{3}&\sqrt{3}&3\\
            1&9&-\sqrt{3}&-3\sqrt{3}&3\\
            12\sqrt{3}&-4\sqrt{3}&0&-8&-4\sqrt{3}\\
            4\sqrt{3}&-12\sqrt{3}&-8&0&4\sqrt{3}\\
            18&18&-6\sqrt{3}&6\sqrt{3}&-2\end{array}\right].
\end{eqnarray*}
Now the equality $f_{Q_{33}^{-1}}=(-1)^j f$ can be written
\begin{eqnarray*}
&&W_5\left[\begin{array}{c}a_1\\
                         a_2\\
                         b_1\\
                         b_2\\
                         c_1\end{array}\right]
=(-1)^j\left[\begin{array}{c}a_1\\
                         a_2\\
                         b_1\\
                         b_2\\
                         c_1\end{array}\right],\
W_2\left[\begin{array}{c}b_5\\
                         b_6\end{array}\right]
=(-1)^j\left[\begin{array}{c}b_3\\
                         b_4\end{array}\right],\\
&&W_2\left[\begin{array}{c}b_3\\
                         b_4\end{array}\right]
=(-1)^j\left[\begin{array}{c}b_5\\
                         b_6\end{array}\right],\
W_3\left[\begin{array}{c}d_1\\
                         d_2\\
                         e\end{array}\right]
=(-1)^j\left[\begin{array}{c}d_1\\
                         d_2\\
                         e\end{array}\right].
\end{eqnarray*}
Suppose $j=1$ and let $\diag[1,1,\sqrt{3},\sqrt{3},1](W_5+E_5)=A\diag[1,1,\sqrt{3},\sqrt{3},3]$. We can easily verify that the the linear span
of the row vectors of $A$ is the linear span of three row vectors $[1,0,0,1/6,5/54]$, $[0,1,0,-1/6,1/54]$ and $[0,0,1,-1,-4/9]$.Thus
\[
 a_1=-\frac{\sqrt{3}}{6}b_2-\frac{5}{18}c_1,\ a_2=\frac{\sqrt{3}}{6}b_2-\frac{1}{18}c_1,\ b_1=b_2+\frac{4\sqrt{3}}{9}c_1,
\]
and $a_1,a_2,b_1,...,b_6,c_1,c_2,d_1,d_2,e$ are linear combinations of $b_2,c_1$. For instance
\begin{eqnarray*}
&&b_3=\sqrt{2}b_2,\ b_4=\sqrt{2}(\sqrt{3}b_2+\frac{8}{9}c_1),\ b_5=\sqrt{2}(\sqrt{3}b_2+\frac{4}{9}c_1),\ b_6=-\sqrt{2}(b_2+\frac{4\sqrt{3}}{9}c_1).
\end{eqnarray*}
Now the condition $W_2[b_4,b_4]=-[b_5,b_6]$ implies that $b_2=c_1=0$, hence $f=0$.

%%Since $W_5+E_5$ is regular, it follows that $a_1=a_2=b_1=b_2=c_1=0$. Hence $b_3=b_4=b_5=b_6=0$, for $f_{Q_{32}^{-1}}=f$.
%% Thus $V(f)$ is singular at $(1,0,0,0)$.
Suppose $j=0$. We can easily show that $\rank(W_5-E_5)=2$ and
$\rank(W_3-E_3)=1$ and that $f_{Q_{31}^{-1}}=f$, $f_{Q_{32}^{-1}}=f$, $f_{Q_{33}^{-1}}=f$ if and only if
\begin{eqnarray*}
&& 2\sqrt{2}a_1-b_5=0,\ 2\sqrt{2}a_2+b_4=0,\ \sqrt{2}b_1+b_6=0,\ 3b_1-d_1=0,\ \sqrt{2}b_2-b_3=0,\\
&& 3b_2-d_2=0,\ c_1=c_2,\ e=-2c_1,\\
&&a_1=\frac{\sqrt{3}}{12}(5b_1+b_2+2\sqrt{3}c_1),\ a_2=\frac{\sqrt{3}}{12 }(-b_1-5b_2+2\sqrt{3}c_1),\ e=\frac{1}{\sqrt{3}}(d_1-d_2),\\
&&b_5=\frac{1}{2}(\sqrt{3}b_3+b_4),\ b_6=\frac{1}{2}(b_3-\sqrt{3}b_4).
\end{eqnarray*}
Consequently $f_{Q_{31}^{-1}}=f$, $f_{Q_{32}^{-1}}=f$, $f_{Q_{33}^{-1}}=f$ if and only if
\begin{eqnarray*}
&&a_1=\frac{\sqrt{3}}{6}(b_1+2b_2),\ a_2=-\frac{\sqrt{3}}{6}(2b_1+b_2),\ b_3=\sqrt{2}b_2,\ b_4=\frac{\sqrt{2}}{\sqrt{3}}(2b_1+b_2),\
b_5=\frac{\sqrt{2}}{\sqrt{3}}(b_1+2b_2),\\
&&b_6=-\sqrt{2}b_1,\ c_1=c_2=\frac{\sqrt{3}}{2}(-b_1+b_2),\ d_1=3b_1,\ d_2=3b_2,\ e=\sqrt{3}(b_1-b_2).
\end{eqnarray*}

\begin{lemma}\label{lem 3.4}
Any $\AAA{5}(4)$-invariant or $\AAA{5}(5)$-invariant quartic surface is singular.
\end{lemma}
\Proof
 Let $f$ be an $\AAA{5}(4)$-invariant quartic form; $f_{Q_{41}^{-1}}\sim f$, $f_{Q_{42}^{-1}}\sim f$ and $f_{Q_{43}^{-1}}\sim f$.
Since $Q_{41}^3=E_4$, $f_{Q_{41}^{-1}}=\omega^i f$ ($i\in [0,2]$). According as $i=0$, $i=1$ or $i=2$, $f$ belongs to
\begin{eqnarray*}
&&\langle x^2z^2,x^2t^2,y^2z^2,y^2t^2,x^2zt,y^2zt,z^2xy,t^2xy,xyzt \rangle,\\
&&\langle x^4,y^4,x^3y,y^3x,z^3x,z^3y,t^3x,t^3y,x^2y^2,z^2xt,z^2yt,t^2xz,t^2yz \rangle,\\
&&\langle z^4,t^4,x^3z,x^3t,y^3z,y^3t,z^3t,t^3z,z^2t^2,x^2yz,x^2yt,y^2xz,y^2xt \rangle,
\end{eqnarray*}
respectively. If $i=0$, then $f$ contains none of monomials $x^4,x^3y,x^3z,x^3t$, hence $V(f)$ is singular at $(1,0,0,0)$.
If $i=1$, $f$ contains none of monomials $x^4,y^4,x^3y,y^3x,z^3x,z^3y,t^3x,t^3y$, for $f_{Q_{43}^{-1}}$ contains
none of them, hence $V(f)$ is singular at $(1,0,0,0,)$. Similar argument shows that $V(f)$ is singular if $i=2$. Similarly
we can show that any $\AAA{5}(5)$-invariant quartic surface is singular.\\

%%%%% section 4
\section{$\SSS{5}$-invariant or $\AAA{6}$-invariant quartic surfaces}

\setcounter{equation}{1}
As is explained in \S2, all faithful representations of $\SSS{5}$ in $PGL_4(k)$ are $\Phi_i$ ($i\in [1,3]$) up to equivalence. They are denoted by
$C_{5!}$\Roman{equation},
\setcounter{equation}{2}
$C_{5!}$\Roman{equation},
\setcounter{equation}{3}
$C_{5!}$\Roman{equation} in \cite{mas}. Let $\SSS{5}(i)=\Phi_i(\SSS{5})$ ($i\in [1,3]$). Similarly, all faithful representations of $\SSS{6}$ in $PGL_4(k)$
are $\varphi_6$ and $\varphi_7$ up to equivalence. Let $\AAA{6}(1)=\varphi_6(\AAA{6})$ and $\AAA{6}(2)=\varphi_7(\AAA{6})$.

 Let
\begin{eqnarray*}
f_0&=& x^4+y^4-2\sqrt{3}x^3y+2\sqrt{3}y^3x+2\sqrt{2}\sqrt{3}z^3x-2\sqrt{2}z^3y+2\sqrt{2}t^3x+2\sqrt{2}\sqrt{3}t^3y\\
    &&+6x^2y^2+6z^2t^2-6\sqrt{3}x^2zt+6\sqrt{3}y^2zt-12xyzt,\\
f_1&=&x^4-y^4+\frac{2}{\sqrt{3}}x^3y+\frac{2}{\sqrt{3}}y^3x+2\frac{\sqrt{2}}{\sqrt{3}}z^3x+2\sqrt{2}z^3y+2\sqrt{2}t^3x\\
    &&-2\frac{\sqrt{2}}{\sqrt{3}}t^3y+2\sqrt{3}(x^2zt+y^2zt).
\end{eqnarray*}
 $V(f_i)$ ($i\in [0,1]$) will be shown to be nonsingular.

\begin{proposition} \label{prop 4.1}
Any  $\SSS{5}(3)$-invariant quartic surface is singular.
An $\SSS{5}(1)$-invariant nonsingular quartic surface  is $V(g_0+\lambda g_1)$.
An $\SSS{5}(2)$-invariant nonsingular quartic surface is $V(f_0)$ or $V(f_1)$.
\end{proposition}
\Proof
 The first part follows from Lemma 3.4, for $\AAA{5}(5)$ is a subgroup of $\SSS{5}(3)$.
By Lemma 3.2 an $\AAA{5}(2)$-invariant quartic form takes the form $g_0+\lambda g_1$, which is $R_{14}$-invariant.
Since  $\SSS{5}(1)=\langle \AAA{5}(2),\ (R_{14})\rangle$, the second part follows.
Let $R_4=iR_{24}$ and $V(f)$ be an $\SSS{5}(2)$-invariant nonsingular quartic surface.
Then $\ord(R_4)=2$,  $\SSS{5}(2)=\langle \AAA{5}(3),\ (R_4)\rangle$, and $f$ is $f_{R_4}=(-1)^j f$ ($j\in [0,1]$). Since $f$
has the form
\begin{eqnarray*}
&&\frac{\sqrt{3}}{6}(b_1+2b_2)x^4-\frac{\sqrt{3}}{6}(2b_1+b_2)y^4\\
&&+b_1x^3y+b_2xy^3+\sqrt{2}b_2xz^3+\frac{\sqrt{2}}{\sqrt{3}}(2b_1+b_2)yz^3+\frac{\sqrt{2}}{\sqrt{3}}(b_1+2b_2)xt^3-\sqrt{2}b_1yt^3\\
&&+\frac{\sqrt{3}}{2}(-b_1+b_2)x^2y^2+\frac{\sqrt{3}}{2}(-b_1+b_2)z^2t^2+3b_1x^2zt+3b_2y^2zt+\sqrt{3}(b_1-b_2)xyzt.
\end{eqnarray*}
by Lemma 3.3, it is easy to see that $f$ is proportional to $f_0$ or $f_1$ according as $j=0$ or $j=1$.\\

\begin{lemma} \label{lem 4.2}
The $V(f_0)$ is nonsingular.
\end{lemma}
\Proof
Denote $f_0$ by $f$, and assume that $f_x$, $f_y$, $f_z$ and $f_t$ vanish at $(a,b,c,d)\in P^3$.
It is easy to see that $cd\not=0$. We may assume $d=1$.
Since $f_x+\sqrt{3}f_y$ and $cf_z-f_t$ vanish, we obtain
\begin{eqnarray*}
&&(\sqrt{3}a-b)c=-\frac{1}{12}a^3+\frac{\sqrt{3}}{4}a^2b+\frac{5}{4}ab^2+\frac{\sqrt{3}}{4}b^3+\frac{\sqrt{2}}{3},\\
&&(\sqrt{3}a-b)c^3=a+\sqrt{3}b.
\end{eqnarray*}
Thus
\[
(\sqrt{3}a-b)f_x=5\sqrt{3}a^4-30a^3b-30ab^3-t\sqrt{3}b^4=5\sqrt{3}(a-ib)(a+ib)(a-\sqrt{3}b-2b)(a-\sqrt{3}b+2b)=0.
\]
Note that $f_z+\sqrt{2}(\sqrt{3}a-b)f_t/2$ takes the form
\[
 c\left[12\sqrt{2}(\sqrt{3}a-b)c +12+3\sqrt{2}(\sqrt{3}a-b)(-3\sqrt{3}a^2-6ab+3\sqrt{3}b^2)\right].
\]
Since $c\not=0$, it follows that $a^3-3ab^2-\sqrt{2}=0$. In particular, $ab\not=0$. Assume $a=ib$ for instance.
Then $c^3=(a+\sqrt{3}b)/(\sqrt{3}a-b)=-i$. The equality $a^3-3ab^2-\sqrt{2}=0$ yields $b^3=i\sqrt{2}/4$, hence
$ab^2=-\sqrt{2}/4$, $a^2b=-\sqrt{2}/4$, and $a^3=\sqrt{2}/4$. Now
\[
 u=-\frac{1}{12}a^3+\frac{\sqrt{3}}{4}a^2b+\frac{5}{4}ab^2+\frac{\sqrt{3}}{4}b^3+\frac{\sqrt{2}}{3}
\]
is equal to zero. Thus $0=u^3=(\sqrt{3}a-b)^3c^3\not=0$, a contradiction. Similarly, $a=-ib$ or
$a=(\sqrt{3}\pm 2)b$ leads to a similar contradiction.\\

\begin{lemma} \label{lem 4.3}
The $V(f_1)$ is nonsingular.
\end{lemma}
\Proof
Denote $f_1$ by $g$, and assume that $g_x$, $g_y$, $g_z$ and $g_t$ vanish at $(a,b,c,d)\in P^3$.
It is easy to see that $bd\not=0$. We may assume $d=1$. Now
\begin{eqnarray*}
&&g_x=4a^3-2\sqrt{3}a^2b+\frac{2}{\sqrt{3}}b^3+2\frac{\sqrt{2}}{\sqrt{3}}c^3+2\sqrt{2}+4\sqrt{3}ac,\\
&&g_y=\frac{2}{\sqrt{3}}a^3+2\sqrt{3}ab^2-4c^3+2\sqrt{2}c^3-2\frac{\sqrt{2}}{\sqrt{3}}+4\sqrt{3}bc,\\
&&g_z=2\sqrt{2}\sqrt{3}(a+\sqrt{3}b)c^2+2\sqrt{2}(a^2+b^2),\\
&&g_t=2\sqrt{2}\sqrt{3}(\sqrt{3}a-b)+2\sqrt{3}(a^2+b^2)c.
\end{eqnarray*}
Since $cg_z-g_t=0$, we obtain $(a+\sqrt{3}b)c^3=\sqrt{3}a-b$. Now the equalities $(a+\sqrt{3}b)g_x-2g_t=0$ and
$(a+\sqrt{3}b)g_y-2\sqrt{3}g_t=0$, and $\sqrt{3}g_x-g_y=0$  can be written
\begin{eqnarray*}
&&\frac{\sqrt{3}}{3}a^4+\frac{1}{2}a^3b-\frac{\sqrt{3}}{2}a^2b^2+\frac{1}{6}ab^3+\frac{\sqrt{3}}{6}b^4-2\frac{\sqrt{2}\sqrt{3}}{3}a+2\frac{\sqrt{2}}{3}b
+b(\sqrt{3}a-b)c=0,\\
&&-\frac{1}{6}a^4-\frac{\sqrt{3}}{6}a^3b-\frac{1}{2}a^2b^2-\frac{\sqrt{3}}{6}ab^3+b^4+\frac{8\sqrt{2}}{3}a-2\frac{\sqrt{2}\sqrt{3}}{3}b
+a(\sqrt{3}a-b)c=0,\\
&&\frac{5}{6}a^3-\frac{\sqrt{3}}{2}a^2b-\frac{1}{2}ab^2+\frac{\sqrt{3}}{2}b^3+2\frac{\sqrt{2}}{3}+(\sqrt{3}a-b)c=0.
\end{eqnarray*}
Substracting the third multiplied by $b$ from the first,  we get
\[
 a^4-\frac{\sqrt{3}}{3}a^3b+2\frac{\sqrt{2}}{\sqrt{3}}ab^3-b^4-2\sqrt{2}a=0.
\]
Substracting the third multiplied by $a$ from the second, we get
\[
a^4-\frac{\sqrt{3}}{3}a^3b+2\frac{\sqrt{2}}{\sqrt{3}}ab^3-b^4-2\sqrt{2}a-2\frac{\sqrt{2}\sqrt{3}}{3}b=0,
\]
hence $b=0$. On the other hand it is easy to see that $(a,0,c,d)\in P^3$ cannot be a singular point of $V(f_1)$.\\

\begin{lemma} \label{lem 4.4}
There exists no $\AAA{6}$-invariant nonsingular quartic surface. Hence, there exists no $\SSS{6}$-invariant nonsingular quartic surface.
\end{lemma}
\Proof
Let $X=Q_{64}$ and $Y=Q_{74}$ (see \S2).
Any subgroup of $PGL_4(k)$ which are isomorphic to $\AAA{6}$ are conjugate to $\AAA{6}(1)$ generated by $\AAA{5}(3)$ and $(X)$ or $\AAA{6}(2)$
generated by $\AAA{5}(5)$ and $(Y)$. Any $\AAA{6}(2)$-invariant quartic surface is singular by Lemma 3.4. Suppose there exists an $\AAA{6}(1)$-invariant
quartic surface $V(f)$. Then $f$ has the form given in Lemma 3.3. But the condition $f_{X^{-1}}\sim f$, namely $f_{X^{-1}}=f$ or $f_{X^{-1}}=-f$,
yields $b_1=b_2=0$.\\

%%% remark 4.5
\begin{remark} Let $h=x^4+y^4+z^4+t^4+12xyzt$, and $G=\Paut(V(h))$. Then $G$ contains a subgroup conjugate to $\SSS{5}(1)$ by Proposition 5.2.
Besides, we can show  by use of computer that there exist no subgroups  conjugate to $\SSS{5}(2)$ in $G$. Consequently none of $V(f_i)$ $(i\in [0,1])$
in Proposition 4.1 is  projectively equivalent to $V(h)$.
\end{remark}

%%%%% section 5
\section{The projective automorphism group of the quartic surface $V(x^4+y^4+z^4+t^4+12xyzt)$}

$\Paut(V(h))$ stands for the projective automorphism group of the nonsingular quartic surface
$V(h)$, where $h=x^4+y^4+z^4+t^4+12xyzt$. Denote this group by $G_{1920}$. It is known that $|G_{1920}|=1920=2^4\times |\SSS{5}|$ \cite[chap.16,\S272]{bur}.
Throughout this section $G$ stands for $G_{1920}$. We will disccuss the relationship between $G$ and $\SSS{5}$. We shall show that $G$ contains a subgroup
conjugate to $\SSS{5}(1)$ and a subgroup conjugate to $\AAA{5}(3)$ and that $G$ contains a normal subgroup $N$ such that $G/N\cong \SSS{5}$.

$G_{1920}$ contains $(C)$, where
\[
 C=\frac{1}{2}\left[\begin{array}{cccc}
                    -i&-i&-1&1\\
                    -i&i &-1&-1\\
                     i& i&-1& 1\\
                     i&-i&-1&-1\end{array}\right]
\]
with $\ord(C)=5$. $G$ contains $(U_i)$ ($i\in [1,3]$) also, where
\begin{eqnarray*}
&&U_1=\frac{1}{2}\left[\begin{array}{cccc}
                       i& 1& i&-1\\
                      -i& 1& i& 1\\
                      -i& 1&-i&-1\\
                       i& 1&-i& 1\end{array}\right],\
U_2=\frac{1}{2}\left[\begin{array}{cccc}
                       1& 1& i&-i\\
                       1&-1&-i&-i\\
                      -i& i& 1& 1\\
                       i& i& 1&-1\end{array}\right],\
U_3=\frac{1}{2}\left[\begin{array}{cccc}
                       0& 0&-i& 0\\
                       0& 0& 0& 1\\
                       i& 0& 0& 0\\
                       0& 1& 0& 0\end{array}\right].
\end{eqnarray*}
Let $S$ be a subset of the group $PGL_4(k)$. We denote by ${PGL_4(k)}_S$ the subgroup
$\{(A)\in PGL_4(k);\ (A)s(A)^{-1}=s\ {\rm for\ any\ }s\in S\}$ of $PGL_4(k)$.

%%%%%%%%

 For a positive integer $n$ and a permutation $\sigma\in \SSS{n}$, $E_n'$ stands for the unit matrix whose $j$-th column will be denoted by $e_j'$, and
$\hat{\sigma}$ the nonsingular matrix $[e_{\sigma(1)}',...,e_{\sigma(n)}']$. Then $\hat{\sigma\tau}=\hat{\sigma}\hat{\tau}$ for any $\tau\in \SSS{n}$,
and the $j$-th column of the matrix
$\hat{\sigma}\diag[a_1,...,a_n]\hat{\tau}$ is equal to $a_{\tau(j)}e_{\sigma\tau(j)}'$, hence
$\hat{\sigma}^{-1}\diag[a_1,...,a_n]\hat{\sigma}=\diag[a_{\sigma(1)},...,a_{\sigma(n)}]$. Consequently
\[
 G_{96}=\{[e_{\sigma(1)},e_{\sigma(2)},e_{\sigma(3)},e_4]\diag[a,b,c,1];\ abc=a^4=b^4=c^4=1,\ \sigma\in \SSS{3}\}.
\]
 is a subgroup of order $96=16|\SSS{3}|$ in $GL_4(k)$ or $PGL_4(k)$ by abuse of notation.
Let $B=\hat{\sigma}$, where $\sigma=(1234)\in \SSS{4}$, and
\begin{eqnarray*}
&&C=\frac{1}{2}\left[\begin{array}{cccc}
                   -i&-i&-1& 1\\
                   -i& i&-1&-1\\
                    i& i&-1& 1\\
                    i&-i&-1&-1\end{array}\right],\
C^2=\frac{1}{2}\left[\begin{array}{cccc}
                    -1&-i& i&-1\\
                    -i&-1& 1&-i\\
                     1&-i&-i&-1\\
                    -i& 1& 1& i\end{array}\right],\\
&&C^3=\frac{1}{2}\left[\begin{array}{cccc}
                    -1& i& 1& i\\
                     i&-1& i& 1\\
                    -i& 1& i& 1\\
                    -1& i&-1&-i\end{array}\right],\
C^4=\frac{1}{2}\left[\begin{array}{cccc}
                     i& i&-i&-i\\
                     i&-i&-i& i\\
                    -1&-1&-1&-1\\
                     1&-1& 1&-1\end{array}\right].
\end{eqnarray*}
We note that $\ord(B)=4$ and  $\ord(C)=5$.

Let $h_\lambda=x^4+y^4+z^4+t^4+\lambda xyzt$. Note that $h_{\lambda}(i^\ell x,y,z)=h_{i^\ell\lambda}(x,y,z)$ for any $\ell\in [0,3]$, namely
$V(h_\lambda)$ and $V(h_{i^\ell\lambda})$ are projectively equivalent.
%%% lemma 5.1
\begin{lemma} \label{lemm 5.1}
$(1)$ The quartic surface $V(h_\lambda)$ is singular if and only if $(\lambda/4)^4=1$.\\
$(2)$ $G_{384}=G_{96}+(B)G_{96}+(B)^2G_{96}+(B)^3G_{96}$ is a group of order $384$.\\
$(3)$ $G_{1920}=G_{384}+(C)G_{384}+(C)^2G_{384}+(C)^3G_{384}+(C)^4G_{384}$ is a group of order $1920$.\\
$(4)$ $\Paut(V(h_{12}))=G_{1920}$, and $\Paut(V(h_\lambda))=G_{384}$ for $\lambda^4\not\in\{0,4^4,12^4\}$.
\end{lemma}
\Proof
It is trivial that $V(h_0)$ is nonsingular. Suppose $\lambda\not=0$ and that $(V(h_\lambda))$ is singular at $(x,y,z,t)$.
Then $4^4xyzt=\lambda^4xyzt\not=0$, hence $\lambda^4=4^4$. Conversely, if $\lambda^4=4^4$, then $V(h_\lambda)$ is singular at
$(1,1,-4/\lambda,1)$.

Let $\mu=\lambda/12$, $g=12^{-4}\Hess(h_\lambda)$, and assume $\mu^4\not\in \{0,3^{-4}\}$.   Then
\[
 g=(1-3\mu^4)x^2y^2z^2t^2+2\mu^3xyzt(x^4+y^4+z^4+t^4)-\mu^2\{x^4(y^4+z^4+t^4)+y^4(z^4+t^4)+z^4t^4\},
\]
so that
\begin{eqnarray*}
&&g_x/2=(1-3\mu^4)xy^2z^2t^2+\mu^3(5x^4yzt+y^5zt+yz^5t+yzt^5)-2\mu^2x^3(y^4+z^4+t^4),\\
 &&g_y/2=(1-3\mu^4)x^2yz^2t^2+\mu^3(x^5zt+5xy^4zt+xz^5t+xzt^5)-2\mu^2y^3(x^4+z^4+t^4),\\
&&g_z/2=(1-3\mu^4)x^2y^2zt^2+\mu^3(x^5yt+xy^5t+5xyz^4t+xyt^5)-2\mu^2z^3(x^4+y^4+t^4),\\
&&g_t/2=(1-3\mu^4)x^2 y^2z^2t+\mu^3(x^5yz+xy^5z+xyz^5+5xyzt^4)-2\mu^2t^3(x^4+y^4+z^4).
\end{eqnarray*}
As is well known, $\Paut(h_\lambda)$ is a subgroup of $\Paut(g)$. Note that $G_{96}\cup (\hat{S_4})\subset \Paut(h_\lambda)$.
Clearly $V(g)$  is singular at $P_1=(1,0,0,0)$, $P_2=(0,1,0,0)$, $P_3=(0,0,1,0)$ and $P_4=(0,0,0,1)$.
If $V(g)$ is singular at $P=(x,y,z,t)$ with $xyzt=0$, then it can be shown easily that $P=P_i$ for some $i\in [1,4]$. Suppose
$V(g)$ is singular at  $P=(x,y,z,t)$ with $xyzt\not=0$. We may assume $t=1$.  Now $g_x=g_y=g_z=g_t=0$ if and only if
$xg_x=yg_y=zg_z=g_t=0$. The condition  $xg_x-g_t=yg_y-g_t=zg_z-g_t=0$ can be written
 \[
 (2\mu xyz-y^4-z^4)(x^4-1)=(2\mu xyz-x^4-z^4)(y^4-1)=(2\mu xyz-x^4-y^4)(z^4-1)=0.
\]
Note that if $x^4=1$, then $y^4=z^4=1$. To see this, first assume $y^4,z^4\not=1$. Then the above equalities imply $y^4=z^4$ and
$2\mu xyz=(\beta+1)$, where $\beta=y^4$. Consequently $16\mu^4 \beta^2=(\beta+1)^4$. Hence the equality $\mu^2 g_t=0$ can be written
$(\beta+1)^2+\mu^4(\beta^2-6\beta+1)=0$, which yields
\[
 0=16\beta^2(\beta+1)^2+(\beta+1)^4(\beta^2-6\beta+1)=(\beta+1)^2(\beta-1)^4,
\]
hence $y^4=\beta=1$, a contradiction. Secondly, assume $y^4=1$ and $z^4\not=1$. Then $2\mu xyz=2$, hence $\mu^4\gamma=1$, where
$\gamma=z^4$. Now the conditon $\mu^2g_t=0$ yields $\gamma=1$, a contardiction. Thus, either $x^4=y^4=z^4=1$ or
$x^4,y^4,z^4\not=1$. We further note that if $x^4=y^4=z^4=1$, then $\mu=i^{-a}$ ($a\in [0,3]$) if $xyz=i^a$. Indeed,
\begin{eqnarray*}
0&=& 12^{-4}g(x,y,z,1)=(1-3\mu^4)i^{2a}+8\mu^3i^a-6\mu^2=-3i^{2a}\{(i^a\mu)^4-\frac{8}{3}(i^a\mu)^3+2(i^a\mu)^2-\frac{1}{3}\}\\
&=&(i^a\mu-1)^3(i^a\mu+\frac{1}{3}),
\end{eqnarray*}
and $(3\mu)^4\not=1$, for $(\lambda/4)^4\not=1$.  Moreover, if $\mu^4=1$, and $x^4=y^4=z^4=1$, then $V(g)$ is singular at $(x,y,z,1)$ if and only if
$\mu xyz=1$, for $\mu xg_x=\mu yg_y=\mu zg_z=\mu g_t=-4(\mu xyz-1)(\mu xyz+3)$, for $\mu xyz+3\not=0$. Next suppose
\[
 2\mu xyz=y^4+z^4,\ \ 2\mu xyz=x^4+z^4,\ \ 2\mu xyz=x^4+y^4.
\]
These three conditions hold if $\mu^4\not=1,\ 0$. Note also that if one of the three conditions does not hold, then
$x^4=y^4=z^4=1$. Now $x^4=y^4=z^4=\mu^4$ so that $x=i^a\mu$, $y=i^b\mu$ and $z=i^c\mu$ such that
$a+b+c=0$ $(\mod 4)$. As far as $\mu^4\not=0,\ 3^{-4}$, the 16 points $(i^a\mu,i^b\mu,i^{-a-b}\mu,1)$ are singular points on $V(g)$.
Indeed, $xg_x$, $yg_y$, $zg_z$ and $g_t$ vanish there. Note also that if $\mu^4=1$, the set of these 16 points in $P^3$ coincides with the set
of 16 points $(x,y,z,1)$ such that $x^4=y^4=z^4=1$ and $\mu xyz=1$. So far, we have shown that
$V(g)$ with $\mu^4\not=0,\ 3^{-4}$ has exactly 20 points $P_i$ ($i\in [1,20]$),
where $\{P_i\ ;\ i\in [5,20]\}=\{(i^a\mu,i^b\mu,i^c\mu,1);\ a+b+c=0\ (\mod\ 4)\}$, points $(i^a\mu,i^b\mu,i^c\mu,1)$ being ordered in dictionary order
such as $P_5=(1,1,1,1)$, $P_6=(1,i,-i,1)$ and $P_{20}=(-i,-i,-1,1)$.
Clearly $\Paut(V(g))$ contains subgroups which act transitively on $\{P_1,P_2,P_3,P_4\}$ or
$\{P_i;\ i\in [5,20]\}$. We will show that if $\mu^4\not=0,1,3^{-4}$, no $(A)\in \Paut(V(g))$ such that $(A)P_4=P_5$, where
$P_4=(0,0,0,1)$ and $P_5=(\mu,\mu,\mu,1)$.  In fact, we can show that the tangent cones \cite[p.79]{sha} to $V(g)$ at $P_4$ and $P_5$ are
$V(g_{P_4})$ and $g_{P_5}$, where $g_{P_4}=xyz$ and $g_{P_5}=\mu^4(\mu^4-1)\{5(x^2+y^2+z^2)-6(xy+yz+zx)\}$. Clearly these two tangent cones,
are affine varieties, are not isomorphic. Thus $\Paut(V(g))$ acts transitively on the set $\{P_1,P_2,P_3,P_4\}$, provided $\mu^4\not=0,1,3^{-4}$.
On the other hand, if $\mu^4=1$, $\Paut(V(g))$ acts transitively on the 20-point set of all singular points of $V(g)$. Indeed, $(C)P_3=P_5$, provided
$\mu=1$. Note that the four projective varieties $V(g)$ with $\mu^4=1$ are projectively equivalent. Indeed, since $V(g)=V(\mu^{-2}g)$, they are
$V(g_{D^{-j}})$ $j\in [0,3]$, where $D=\diag[i,1,1,1]$.

Let $H_{P_4}=\{(A)\in \Paut(V(h_\lambda))\ :\ (A)P_4=P_4\}$ for $\lambda$ such that
$\lambda^4\not\in \{0,4^4,12^4\}$ or $\lambda=12$. $\Paut(V(h_\lambda))$ acts transitively on $\{P_j\ :\ j\in [1,4]\}$ or
$\{P_j\ :\ j\in [1,20]\}$ according as $\lambda^4\not\in\{0,4^4,12^4\}$ or $\lambda=12$, for $\Paut(V(h_\lambda))$ contains $G_{384}$ and
$\Paut(V(h_{12}))$ contains $G_{1920}$. It remains to show $H_{P_4}=G_{96}$, but it suffices to show $H_{P_4}\subset G_{96}$.
Assume $(A)\in H_{P_4}$, where $A=[a_{ij}]\in GL_4(k)$ with $a_{14}=a_{24}=a_{34}=a_{44}-1=0$. As noted in the preliminaries $[0,0,0,1]A\sim
[0,0,0,1]$, hence $a_{41}=a_{42}=a_{43}=0$. Now the condition $h_{\lambda,A^{-1}}\sim h_\lambda$ yields
$A=[e_{\sigma(1)},e_{\sigma(2)},e_{\sigma(3)},e_4]\diag[a,b,c,1]$, where $\sigma\in \SSS{3}$ and $a,b,c\in k^*$, hence
$a^4=b^4=c^4=abc=1$. Thus $(A)\in G_{P_4}$, as desired.\\

  Let $H_1=\diag[1,1,1,1]$, $H_2=\diag[-1,-1,1,1]$, $H_3=\diag[-1,1,-1,1]$, $H_4=\diag[1,-1,-1,1]$, and let $K_1=H_1$,
\begin{eqnarray*}
&&K_2=\left[\begin{array}{cccc}0&1&0&0\\
                               1&0&0&0\\
                               0&0&0&1\\
                               0&0&1&0\end{array}\right],\
K_3=\left[\begin{array}{cccc}0&0&1&0\\
                             0&0&0&1\\
                             1&0&0&0\\
                             0&1&0&0\end{array}\right],\
K_4=\left[\begin{array}{cccc}0&0&0&1\\
                             0&0&1&0\\
                             0&1&0&0\\
                             1&0&0&0\end{array}\right].
\end{eqnarray*} One can easily see that $(H_i)$ and $(K_j)$ commute.
Since $(H_i)$, and  $(K_i)$ ($i\in [1,4]$) form Klein's fourgroups, and $(H_iK_j)$ ($i,j\in [1,4]$) are distinct,
the 16 transformations $(H_iK_j)$ form an abelian subgroup ${\cal{A}}_{16}$ of $G=\Paut(V(h_{12}))$. Let $A_{i+4(j-1)}=H_iK_j$.

In view of Lemma 5.1 we can search for subgroups of $G_{1920}$ which are isomorphic to $\SSS{5}$ or $\AAA{5}$ using computer.
In the proof of Lemma 2.2 (1) we defined the representation $\Psi_1$ of $\SSS{5}$ in $GL_4(k)$ such that $\Psi_(s_j)=R_{ij}$ ($j\in [1,3]$) and
$\Psi_1(t_1)=R_{14}$, hence $\Psi_1(t_j)=r_j$ ($j\in [1,4]$), where $r_1=R_14=[e_1,e_2,e_4,e_3]$, $r_2=[e_1,e_2,\omega^2e_4,\omega e_3]$ and
\begin{eqnarray*}
&&r_3=\left[\begin{array}{cccc}1&0&0&\\
                               0&-\frac{1}{3}&\frac{2}{3}&\frac{2}{3}\\
                               0&\frac{2}{3}&\frac{2}{3}&-\frac{1}{3}\\
                               0&\frac{2}{3}&-\frac{1}{3}&\frac{2}{3}\end{array}\right],\
r_4=\left[\begin{array}{cccc}-\frac{1}{4}&\frac{\sqrt{15}}{4}&0&0\\
                             \frac{\sqrt{15}}{4}&\frac{1}{4}&0&0\\
                             0&0&0&1\\
                             0&0&1&0\end{array}\right].
\end{eqnarray*}
Recall that the representation $\Psi_1$ is equivalent to the four-dimensional irreducible representation $V$ \cite[p.28]{ful}.
Let
\[
 T=\diag[-\frac{\sqrt{15}}{20},\frac{1}{4},1,1]\left[\begin{array}{cccc}
                                               3+i     &1-3i      &3+i &7-i\\
                                               1+3i    &3-i       &1+3i&-3-3i\\
                                               \omega^2&-i\omega  &1   &0\\
                                               \omega  &-i\omega^2&1   &0
                                                                         \end{array}\right].
\]
Note that $T\in GL_4(k)$.
%% proposition 5.2
\begin{proposition}
$(1)$ $G_{1920}$ contains a subgroup  $(T^{-1})\SSS{5}(1)(T)$.\\
$(2)$ $\alpha h_{T}=g_0+\lambda g_1$, where $\alpha=27(-1+3i)/2$ and $\lambda=(-8+9i)/4$.
\end{proposition}
\Proof
Define matrices $F_j\in GL_4(k)$ ($j\in [1,4]$) as follows; $F_1=[ie_2,-ie_1,e_3,e_4],\ F_2=[e_3,e_2,e_1,e_4],$ and
\begin{eqnarray*}
&&\ F_3=\frac{1}{2}\left[\begin{array}{cccc}1&i&-i&1\\
                                          -i&1&1&i\\
                                           i&1&1&-i\\
                                           1&-i&i&1\end{array}\right],
\ F_4=\frac{1}{2}\left[\begin{array}{cccc}1&i&-1&i\\
                                          -i&1&-i&-1\\
                                          -1&i&1&i\\
                                          -i&-1&-i&1\end{array}\right].
\end{eqnarray*}
(1) It can be easily seen by Theorem 1.1 that there exists a faithful representation $\Psi$ of $\SSS{5}$ in $GL_4(k)$ such that
$\Psi(t_j)=F_j$. The character of $\Psi$ coincides with that of the four-dimensional irreducible representation $V$ \cite[p.28]{ful}. Namely
the representations $\Psi$ and $V$ are equivalent. Setting  $R_1=F_1F_2$, $R_2=F_1F_3$, and $R_3=F_1F_4$, we have $\Psi(s_j)=R_j$ ($j\in [1,3]$).
Clearly $\Phi=\pi\circ\Psi$ is a faithful representation of $\SSS{5}$ in $PGL_4(k)$ equivalent to $\Phi_1$. On the other hand, we can verify
that $h_{{F_j}^{-1}}=h$ for $j\in [1,3]$. Since $F_4=\diag[i,-i,i,-i]F_4'[e_2,e_1,e_3,e_4]$, where $F_4'$ is equal to $F_3$ up to the order
of row vectors,
$h_{{F_4}^{-1}}=h$ also holds. Hence $\Phi(\SSS{5})$ is a subgroup of $G_{1920}$, and it is conjugate to $\SSS{5}(1)=\Phi_1(\SSS{5})$.
Since the representations $\Psi_1$ and $\Psi$ are equivalent, there exists an $S=[s_{ij}]\in GL_4(k)$ such that $\Psi_1(\sigma)S=S\Psi(\sigma)$ for
any $\sigma\in \SSS{5}$. The conditions $r_jS=SF_j$ for $j=1,2,4,3$ imply that $S=\omega s_{31}T$ ($s_{31}\in k^*$). Consequently $G_{1920}$
contains $\pi\circ\Psi(\SSS{5})=(T^{-1})(\SSS{5}(1)(T)$.
(2) Since $V(h_T)$ is nonsingular $\SSS{5}(1)$-invariant quarticsurface, $h_T=\alpha^{-1}( g_0+\lambda g_1)$ for some $\alpha,\lambda\in k$, that is,
$\alpha h=(g_0+\lambda g_1)_{T^{-1}}$, by Lemma 3.2.  The left-hand side is equal to the sum of the following six polynomials $p_j$ ($j\in [1,6]$):
\begin{eqnarray*}
&&(-x^3+2\sqrt{15}y^3)x_{T^{-1}}+\{(z+t)^3-3(z+t)zt\}(2\sqrt{15x+10y})_{T^{-1}}+13x^2y^2_{T^{-1}}\\
&&+5z^2t^2_{T^{-1}}+(26x^2-6\sqrt{15}+20y^2)zt_{T^{-1}}+\lambda (x^2+y^2+2zt)^2_{T^{-1}}.
\end{eqnarray*}
We can easily obtain the coefficients of $x^2y^2$ for the $p_j$, whose sum must vanish. Thus
\[
 \frac{-13122-16704i}{1600}+12(-1+3i)+\frac{585}{32}-15+0+\frac{36(-7+i)}{25}\lambda=0,
\]
which gives $\lambda=(-8+9i)/4$. We can also obtain the coefficeints $c_j$ of $t^4$ for $p_j$ as follows.
\[
 c_1=\frac{-16893+19224i}{1600},\ c_2=0,\ c_3=\frac{2457+8424i}{320},\ c_4=0,\ c_5=0,\ c_6=\frac{-8+9i}{4}\cdot \frac{72+54i}{25}.
\]
Thus $\alpha=\sum_{j=1}^6c_j=27(-1+3i)/2$.\\

%% proposition 5.3
\begin{proposition}
$G_{1920}$ contains a subgroup conjugate to $\AAA{5}(3)$.
\end{proposition}
\Proof
Define matrices $Q_j\in GL_4(k)$ ($j\in [1,3]$) as follows.
\[
 Q_1=-\frac{1}{2}\left[\begin{array}{cccc}
-1&i&1&i\\
1&i&1&-i\\
1&i&-1&i\\
-1&i&-1&-i\end{array}\right],\
Q_2=-\frac{1}{2}\left[\begin{array}{cccc}
-1&-1&1&1\\
-1&1&1&-1\\
1&1&1&1\\
1&-1&1&-1\end{array}\right],\
Q_3=-\frac{1}{2}\left[\begin{array}{cccc}
-1&-1&1&1\\
-1&1&-1&1\\
1&-1&-1&1\\
1&1&1&1\end{array}\right].
\]
By Theorem 1.1 there exists a faithful representation $\varphi$ of $\AAA{5}$ in $PGL_4(k)$ such that
$\varphi(s_j)=Q_j$ ($j\in [1,3]$). We shall show that $\varphi$ and $\varphi_3$ defined in \S2 are
equivalent. It suffices to show that $S^{-1}Q_jS=Q_{3j}$ ($j\in [1,3]$) for some $S\in GL_4(k)$.
Let
\begin{eqnarray*}
&&T=\left[\begin{array}{cccc}
1&1-i&i\omega-\omega^2&i\omega^2-\omega\\
0&i&-2i+\omega-\omega^2&-2i+\omega^2-\omega\\
-1&0&i\omega-\omega^2&i\omega^2-\omega\\
0&1&1&1\end{array}\right],\\
&&U=\left[\begin{array}{cccc}
1+\sqrt{3}+i(2+\sqrt{3})&1-\sqrt{3}+i(2-\sqrt{3})&0&0\\
1&                       1&                       0&0\\
0&                       0&                       1&0\\
0&                       0&                       0&1\end{array}
\right],
\end{eqnarray*}
and $V=\diag[-\sqrt{2}\omega^2,\sqrt{2}(2+\sqrt{3})\omega^2,(2+\sqrt{3})\omega,1]$, and note that
\begin{eqnarray*}
&&T^{-1}=\frac{1}{12}\left[\begin{array}{cccc}
4&2(1+i)&-8&2(-1+i)\\
2(1+i)&-4i&2(1+i)&4\\
(-\sqrt{3}-1)(1+i)&2&(-\sqrt{3}-1)(1+i)&4+2\sqrt{3}\\
(\sqrt{3}-1)(1+i)&2&(\sqrt{3}-1)(1+i)&4-2\sqrt{3}\end{array}
\right],\\
&&U^{-1}= \frac{1}{4\sqrt{3}}\left[\begin{array}{cccc}
1-i&-3+2\sqrt{3}-i& 0& 0\\
-1+i&3+2\sqrt{3}+i& 0& 0\\
0& 0&             4\sqrt{3}&0\\
0& 0&               0&4\sqrt{3}\end{array}
\right],
\end{eqnarray*}
and $V^{-1}=\diag[-\frac{\omega}{\sqrt{2}},\frac{(2-\sqrt{3})\omega}{\sqrt{2}},(2-\sqrt{3})\omega^2,1]$.
Now $T^{-1}Q_1T=Q_{31}$. Moreover, denoting the $i$-th row of $12T^{-1}Q_{2}T$ by $q_i$ we get \\
\ \ \ \ $q_1=[2(3+i),2(5-i),6-4\sqrt{3}+i(-12+10\sqrt{3}),6+4\sqrt{3}+i(-12-10\sqrt{3})]$,\\
\ \ \ \ $q_2=[2(1-i),-2(3+i),2\sqrt{3}-6i,-2\sqrt{3}-6i]$,\\
\ \ \ \ $q_3=[-1-\sqrt{3}+i(1-\sqrt{3}),3+\sqrt{3}+i(7+3\sqrt{3}),4\sqrt{3},0]$,\\
\ \ \ \ $q_4=[-1+\sqrt{3}+i(1+\sqrt{3}),3-\sqrt{3}+i(7-3\sqrt{3}),0,-4\sqrt{3}]$.\\
$T^{-1}Q_{3}T$ takes the form
\begin{eqnarray*}
%&&\frac{1}{12}\left[\begin{array}{cccc}
%2(3+i)&2(5-i)&6-4\sqrt{3}+i(-12+10\sqrt{3})&6+4\sqrt{3}+i(-12-10\sqrt{3})\\
%2(1-i)&-2(3+i)&2\sqrt{3}-6i&-2\sqrt{3}-6i\\
%-1-\sqrt{3}+i(1-\sqrt{3})&3+\sqrt{3}+i(7+3\sqrt{3})&4\sqrt{3}&0\\
%-1+\sqrt{3}+i(1+\sqrt{3})&3-\sqrt{3}+i(7-3\sqrt{3})&0&-4\sqrt{3}\end{array}
%\right],\\
%&&\frac{1}{12}\left[\begin{array}{cccc}
%2(3+i)&2(5-i)&4\{\omega-4\omega^2-2i(1-\omega)\}&4\{\omega^2-4\omega-2i(1-\omega^2)\}\\
%2(1-i)&-2(3+i)&4(\omega^2-1)i&4(\omega-1)i\\
%-1-\sqrt{3}+i(1-\sqrt{3})&3+\sqrt{3}+i(7+3\sqrt{3})&4\sqrt{3}&0\\
%-1+\sqrt{3}+i(1+\sqrt{3})&3-\sqrt{3}+i(7-3\sqrt{3})&0&-4\sqrt{3}\end{array}
%\right], \\
&&\frac{1}{2}\left[\begin{array}{cccc}
2&2(1-i)&0&0\\
0&-2&0&0\\
0&0&0&(2+\sqrt{3})(-1+\sqrt{3}i)\\
0&0&(-2+\sqrt{3})(1+\sqrt{3}i)&0\end{array}
\right].
\end{eqnarray*}
Therefore, $S^{-1}Q_jS=Q_{3j}$ ($j\in [1,3]$) for $S=TUV$.\\
To see that $\varphi(\AAA{5})$ lies in $G_{1920}$, we note that
\begin{eqnarray*}\diag[i,-1,i,1]Q_1[e_1,e_2,e_4,e_3],\ \ \ {\rm and\ \ \  }
\diag[-i,i,1,1]Q_3\diag[1,i,-i,1]
\end{eqnarray*}
are equal to $R_3$ up to the order of row vectors and
that $Q_2=[e_1,e_2,e_4,e_3]Q_3[e_1,e_2,e_4,e_3]$. Consequently $h_{{Q_j}^{-1}}=h$ ($j\in [1,3]$). \\

%% proposition 5.4
\begin{proposition}\label {pro 5.4} \  \\
$(1)$ ${\cal{A}}_{16}\triangleleft G_{1920}$.\\
$(2)$ $G_{1920}/{\cal{A}}_{16}\cong \SSS{5}$.
\end{proposition}
\Proof
(1) It suffices to show that $(T)(H_j)(T)^{-1}$ and $(T)(K_j)(T)^{-1}$ ($j\in [2,4]$) belong to ${\cal{A}}_{16}$
for $T=\diag[i^a,i^b,i^{-c},1]$ ($a+b=-c\ (\mod 4\ )$), $T=\hat{\sigma}$ ($\sigma\in \SSS{3}\subset \SSS{4}$), $T=B$ and $T=C$. For example,
if $T=\diag[i^a,i^b,i^c,1]$, then
\begin{eqnarray*}
&&TK_2T^{-1}=\left[\begin{array}{cccc}0&i^{a-b}&0&0\\
                                      i^{b-a}&0&0&0\\
                                      0&0&0&i^{c}\\
                                      0&0&i^{-c}&0\end{array}\right]
\sim \left[\begin{array}{cccc}        0&i^{a-b+c}&0&0\\
                                      i^{b-a+c}&0&0&0\\
                                      0&0&0&i^{2c}\\
                                      0&0&1&0\end{array}\right]\\
&&= \left[\begin{array}{cccc}        0&i^{2b}&0&0\\
                                      i^{2a}&0&0&0\\
                                      0&0&0&i^{2(a+b)}\\
                                      0&0&1&0\end{array}\right]
= \left[\begin{array}{cccc}        0&(-1)^{b}&0&0\\
                                      (-1)^{a}&0&0&0\\
                                      0&0&0&(-1)^{(a+b)}\\
                                      0&0&1&0\end{array}\right].
\end{eqnarray*}
Thus $(T)(K_2)(T)^{-1}$ belongs to ${\cal{A}}_{16}$.
(2)  Since ${\cal{A}}_{16}P_1=\{P_j\ :\ j\in [1,4]\}$, ${\cal{A}}_{16}\triangleleft G$, and
$G=G_{1920}$ acts transitively on ${\cal P}=\{P_i;\ i\in [1,20]\}$ (the set of all singular points on $V(\Hess(h_{12}))$),
 ${\cal P}$ is the union of five ${\cal{A}}_{16}$-orbits and we have a group homomorphism $\varphi$
from $G$ to the permutation group of these orbits. To be more prcise, ${\cal{P}}={\cal{P}}_1+{\cal{P}}_2+{\cal{P}}_3+{\cal{P}}_4+{\cal{P}}_5$,
where
\begin{eqnarray*}
 {\cal{P}}_1&=&\{P_1,P_2,P_3,P_4\},\  {\cal{P}}_2=\{P_5,P_7,P_{13},P_{15}\},\  {\cal{P}}_3=\{P_9,P_{11},P_{17},P_{19}\},\\
{\cal{P}}_4&=&\{P_6,P_8,P_{14},P_{16}\},\  {\cal{P}}_5=\{P_{10},P_{12},P_{18},P_{20}\},
\end{eqnarray*}
%%%%%%% original order of {\cal{P}}_i
%\begin{eqnarray*}
% {\cal{P}}_1&=&\{P_1,P_2,P_3,P_4\},\  {\cal{P}}_2=\{P_5,P_7,P_{13},P_{15}\},\  {\cal{P}}_3=\{P_6,P_8,P_{14},P_{16}\},\\
%{\cal{P}}_4&=&\{P_9,P_{11},P_{17},P_{19}\},\  {\cal{P}}_5=\{P_{10},P_{12},P_{18},P_{20}\},
%\end{eqnarray*}
and $(T){\cal{P}}_i={\cal{P}}_{\sigma(i)}$, where $\sigma=\varphi((T))\in \SSS{5}$. It remains to show that $\Ker(\varphi)={\cal{A}}_{16}$, for $|G|=1920$.
Let $\sigma=\varphi((C))$ and $\tau=\varphi(([e_2,e_1,e_3,e_4]))$. We can easily verify that $\sigma=(12345)$ and $\tau=(34)$. Thus $\Ima\ \varphi$ contains
transpositions $(j\ j+1)$ ($j\in [1,4]$), hence $\varphi$ is surjective so that $|\Ker(\varphi)|=|G_{1920}|/|\SSS{5}|=16$.\\

\end{document}